\documentclass[11pt]{article}
\usepackage[a4paper,margin=2.5cm]{geometry}
\usepackage{amsmath,amssymb,amsfonts,amsthm,mathtools}
\usepackage{bbm}
\usepackage{enumitem}
\usepackage{hyperref}
\usepackage{xcolor}
\usepackage{float}

\usepackage{hyperref}
\hypersetup{
    colorlinks,
    citecolor=blue,
    linkcolor=red,
    urlcolor=blue
}

\newtheorem{theorem}{Theorem}
\newtheorem{prop}{Proposition}

\newtheorem{cor}{Corollary}

\newtheorem{remark}{Remark}

\DeclareMathOperator{\MMD}{MMD}

\title{A Kernel Two-Sample Test Invariant under Group Action with Applications to Functional Data}

\author{
Madison Giacofci\thanks{Université Rennes 2, IRMAR, UMR CNRS 6625, Rennes, France}
\and 
Anouar Meynaoui\footnotemark[1]
\and
Alex Podgorny\thanks{Ensai, CNRS, CREST-UMR 9194, Rennes, France} 
}

\date{}

\begin{document}

\maketitle

\begin{abstract}
We introduce a kernel-based two-sample test for comparing probability distributions up to group actions. Our construction yields invariant kernels for locally compact $\sigma$-compact groups and extends classical Haar-based approaches beyond the compact setting. The resulting invariant Maximum Mean Discrepancy (MMD) test is developed in a general framework where the sample space is assumed to be Polish. Under natural conditions, the invariant kernel induces a characteristic kernel on the quotient space, ensuring consistency of the associated MMD test. The method is well suited to functional data, where invariances such as temporal shifts arise naturally, and its effectiveness is illustrated through simulation studies.
\end{abstract}

\section{Introduction}

In many real-world applications, the variability observed in data is partially explained by nuisance transformations. For functional data, typical examples include translations, rotations, scalings, time-shifts, or more general reparametrizations that preserve the underlying content of a given observation. Such transformations arise naturally in many application domains. For instance, images representing the same object or scene may differ by rotations, translations, or changes in scale depending on the camera viewpoint. Similarly, handwritten digits may involve small deformations while still representing the same digit \cite{lecun2002gradient}. In growth curve analysis, individuals may experience biological events such as growth spurts at different ages, producing curves that are essentially identical up to a temporal reparametrization \cite{ramsay2005functional}. In audio analysis, two signals corresponding to the same sound may differ only by a time shift or by a change in duration due to recording conditions \cite{sakoe2003dynamic}. In biomedical signal analysis, electrocardiogram (ECG) signals record the electrical activity of the heart over time. Two ECG signals may look different simply because the heart beats slightly faster or slower, which shifts the timing of the main peaks in the signal. In this case, the overall pattern remains the same, but it appears stretched or shifted in time \cite{shorten2014use}. In such settings, directly comparing the distributions of two datasets can be misleading, since apparent differences may arise solely from nuisance transformations, even when the two samples represent the same underlying phenomenon. This motivates the development of statistical procedures that are invariant to prescribed transformations, so that only meaningful differences between distributions are detected. In this work, we focus specifically on the two-sample testing problem.

\medskip
Kernel-based two-sample tests provide a powerful and flexible framework for comparing probability distributions. In particular, methods based on the maximum mean discrepancy (MMD) have become widely used due to their strong theoretical guarantees and their ability to handle complex and high-dimensional data \cite{gretton2012kernel}. However, standard kernel two-sample tests are typically sensitive to transformations of the data, even though in many applications distributions should be regarded as equal up to nuisance transformations. A natural way to address this issue is to incorporate invariance into the comparison procedure through the action of a group on the observation space. Under this perspective, observations that differ only by such transformations are regarded as equivalent, and the relevant object becomes the distribution induced on the corresponding quotient space. The resulting testing problem is therefore to determine whether two distributions remain different once these transformations are disregarded. A classical way to enforce such invariance in kernel methods consists in averaging a base kernel along the transformations of the group. When the group is compact, this averaging can be performed using the Haar probability measure, leading to kernels that are invariant under the prescribed transformations. This idea has been explored in the machine-learning literature \cite{Haasdonk2005, Mroueh2015, Raj2017}. The same averaging mechanism also appears prominently in the theory of data augmentation. In practice, augmentation replaces each observation with randomly transformed versions. In many pipelines, training proceeds by repeatedly sampling such transformations. Augmenting inputs and then learning with a kernel method is closely related to replacing the original kernel with an augmentation-averaged version, obtained by averaging the base kernel over all pairs of transformations of the two inputs. This relationship is made explicit in several works \cite{Mroueh2015,Dao2019,Chen2020}. These connections suggest that augmentation can be viewed as transforming the underlying distributions before comparison. Another related but distinct problem is to test whether the underlying distribution of a sample is invariant under the action of a given group. The recent work of \cite{Soleymani2025} proposes kernel-based tests for this problem when the acting group is compact. Our objective is different, we instead compare two distributions modulo a group action, i.e., we test equality on the quotient space. 

\medskip
As mentioned earlier, existing approaches rely on the compactness of the transformation group. In many situations of practical interest, transformations such as translations or scalings involve non-compact groups. For such groups, the Haar measure is not finite and cannot be normalized into a probability measure. As a consequence, the averaging construction described above cannot be applied directly, and extending kernel-based testing procedures to such settings requires different ideas. 

\medskip
\textbf{Contributions.}  Our main contributions are as follows.
\begin{itemize}
  \item \textbf{Invariant kernels beyond compact groups.}
  We introduce a weighted averaging procedure for locally compact $\sigma$-compact groups, which yields well-defined invariant kernels. This extends the classical Haar-integration kernel construction, which is limited to compact groups \cite{Haasdonk2005,Mroueh2015}. In the non-compact setting, it is typically replaced by quasi-invariant surrogates \cite{Raj2017}.
  \item \textbf{A rigorous invariant MMD two-sample test.}
  Using these invariant kernels, we formalize an MMD-based invariant two-sample test. Then, we study its statistical properties in a general theoretical framework. In particular, the only assumption on the data space is that it is Polish.
  \item \textbf{Characteristic kernels on the quotient.}
  We show that the invariant kernel induces a kernel on the quotient space. Under natural conditions, this kernel is characteristic. The associated MMD is therefore zero if and only if the two distributions are equal on the quotient. Consequently, the resulting two-sample test is consistent.
\end{itemize}

\medskip
\textbf{Organization of the paper.} Section \ref{sec:background} recalls some background material. In Section \ref{sec:kernel_review}, we review basic notions on RKHS and kernels, as well as the associated nonparametric two-sample tests. In Section \ref{sec:group-actions}, we recall basic notions on groups, group actions, and the Haar measure. Section \ref{sec:invariant-mmd} presents our invariant MMD two-sample test. It contains the main theoretical contributions of the paper and discusses practical aspects of implementing the test. Section \ref{sec:simulations} presents simulation studies on synthetic signals with temporal shifts. Section \ref{sec:real_data} illustrates the methodology on a real-data application to phonocardiogram (PCG) signals. 

\section{Background}
\label{sec:background}

\subsection{Review on Kernel two-sample tests}
\label{sec:kernel_review}

Nonparametric two-sample testing is a fundamental problem in statistics, where the aim is to determine whether two samples are drawn from the same underlying probability distribution. Mathematically, assume we have two independent i.i.d. samples
$$\{X_i\}_{i=1}^n \sim P, \quad \{Y_i\}_{i=1}^m \sim Q,$$
where $P$ and $Q$ are probability measures on a measurable space $(\mathcal{X},\mathcal{A})$. The two-sample problem consists in testing 
$$\mathcal{H}_0 : P = Q \quad \mbox{against} \quad \mathcal{H}_1 : P \neq Q.$$

Early nonparametric two-sample tests include the Kolmogorov--Smirnov one \cite{hodges1958significance}, which compares empirical cumulative distribution functions, via the supremum norm. A major limitation of this test is that it is restricted to one-dimensional data, as its formulation relies on the existence of a natural total ordering of the observations. Other classical two-sample procedures are based on the Cramer--von Mises criterion, which replaces the supremum by an integrated squared difference \cite{anderson1962cvm}. A closely related test is the Anderson--Darling procedure, originally developed in the one-sample setting \cite{andersondarling1952,andersondarling1954} and later extended to the two-sample problem \cite{pettitt1976two}, where a weighting function is introduced in the integrated squared difference. Like the Kolmogorov-Smirnov test, these procedures fundamentally exploit univariate ordering and therefore do not admit a canonical multivariate extension. To overcome these limitations, kernel-based two-sample tests have been introduced. These approaches compare probability distributions by mapping them into a Reproducing Kernel Hilbert Space (RKHS), thereby reducing the two-sample test to the comparison of elements in a Hilbert space. The discrepancy induced by this embedding leads to the Maximum Mean Discrepancy (MMD), which defines a metric on the space of probability measures for the so-called characteristic kernels \cite{gretton2006kernel,gretton2012kernel}. Kernel-based two-sample tests offer several advantages, including the ability to detect general distributional differences in arbitrary dimensions, theoretical guarantees of consistency for characteristic kernels, and flexibility through kernel choice \cite{gretton2012kernel}. Moreover, the MMD admits simple empirical estimators with good statistical properties and can be efficiently computed in practice.  Accordingly, we focus on MMD-based testing procedures in the remainder of this work.

\medskip
Before introducing kernel-based two-sample tests, we first recall some key notions for their construction. Let $\mathcal{X}$ be a non-empty set and let $\mathcal{H}$ be a Hilbert space of real-valued functions on $\mathcal{X}$ equipped with inner product $\langle\cdot,\cdot\rangle_{\mathcal{H}}$ and norm $\|\cdot\|_{\mathcal{H}}$. The space $\mathcal{H}$ is called a Reproducing Kernel Hilbert Space (RKHS) if, for every $x\in\mathcal{X}$, the point-evaluation functional $L_x:\mathcal{H}\to\mathbb{R}$ defined by $L_x(f)=f(x)$ is continuous. By the Riesz representation theorem, this implies that for each $x\in\mathcal{X}$ there exists a unique element $\phi(x)\in\mathcal{H}$ such that $f(x)=\langle f,\phi(x)\rangle_{\mathcal{H}}$ for all $f\in\mathcal{H}$, known as the reproducing property \cite{aronszajn1950theory}. The associated reproducing kernel is defined as $k:\mathcal{X}\times\mathcal{X}\to\mathbb{R}$, $k(x,x'):=\langle \phi(x),\phi(x')\rangle_{\mathcal{H}}$, and satisfies $k(\cdot,x)\in\mathcal{H}$ together with $f(x)=\langle f,k(\cdot,x)\rangle_{\mathcal{H}}$ for all $f\in\mathcal{H}$. An RKHS is uniquely characterized by its reproducing kernel, and conversely every positive definite kernel defines a unique RKHS. Consequently, we refer to these two notions interchangeably in the following.

\medskip
Beyond representing individual points $x\in\mathcal{X}$ in a RKHS, the kernel framework also allows probability measures on $\mathcal{X}$ to be embedded as elements of the Hilbert space via expectations of feature maps \cite{smola2007hilbert,sriperumbudur2010hilbert}. Given a probability measure $P$ on a measurable space $(\mathcal{X},\mathcal{A})$ such that $\mathbb{E}_{X\sim P} [\sqrt{k(X,X)}]<\infty$, its kernel mean embedding is defined as $m(P) := \mathbb{E}_{X\sim P}[k(\cdot,X)]\in\mathcal{H}$. A kernel $k$ is said to be characteristic if the mean embedding map $P\mapsto m(P)$ is injective, that is, $m(P)=m(Q)$ if and only if $P=Q$. Characteristic kernels therefore ensure that probability distributions are uniquely represented by their mean embeddings in the RKHS, a property that is crucial for defining distances between distributions based on kernel embeddings \cite{sriperumbudur2010hilbert}. Given two probability measures $P$ and $Q$ on $(\mathcal{X},\mathcal{A})$, the Maximum Mean Discrepancy (MMD) associated with a kernel $k$ is defined as $\MMD_k(P,Q):=\|m(P)-m(Q)\|_{\mathcal{H}}$. When the kernel $k$ is characteristic, this quantity defines a metric on the space of probability measures, in the sense that $\MMD_k(P,Q)=0$ if and only if $P=Q$ \cite{gretton2012kernel}. The two-sample problem reduces to testing
$$\mathcal{H}_0: \MMD_k(P,Q)=0 \quad \mbox{against} \quad \mathcal{H}_1:\MMD_k(P,Q)>0.$$
In practice, the $\mathrm{MMD}_k(P,Q)$ is unknown and is estimated from the available samples. Indeed, the MMD admits a closed-form expression in terms of expectations, namely 
$$\MMD_k^2(P,Q) = \mathbb{E}_{X,X'\sim P}[k(X,X')] + \mathbb{E}_{Y,Y'\sim Q}[k(Y,Y')] - 2\mathbb{E}_{X\sim P,Y\sim Q}[k(X,Y)],$$
where all expectations are taken over independent copies. Replacing these expectations by their empirical counterparts yields an unbiased estimator of $\MMD_k^2(P,Q)$ in the form of a two-sample U-statistic of order two. This estimator is given by
$$\widehat{\MMD}_k^2=\frac{1}{n(n-1)}\sum_{i\neq i'}k(X_i,X_{i'})+\frac{1}{m(m-1)}\sum_{j\neq j'}k(Y_j,Y_{j'})-\frac{2}{nm}\sum_{i=1}^n\sum_{j=1}^m k(X_i,Y_j).$$
To construct a level-$\alpha$ test, a natural choice is to reject $\mathcal{H}_0$ for large values of $\widehat{\MMD}_k^2$ by comparing it to the $(1-\alpha)$-quantile of its distribution under $\mathcal{H}_0$. In practice, this null quantile is unknown because the distributions $P$ and $Q$ are unknown, and it is therefore estimated either by a permutation procedure or by approximating the asymptotic distribution under $\mathcal{H}_0$. These different procedures are described below. 

\medskip
The permutation procedure relies on the fact that the pooled sample $(X_1,\dots,X_n,Y_1,\dots,Y_m)$ is exchangeable under $\mathcal{H}_0$. Let $\mathcal{Z}=(Z_1,\dots,Z_{n+m})$ denote this pooled sample with $Z_i=X_i$ for $i\leq n$ and $Z_{n+j}=Y_j$ for $j\leq m$. Under $\mathcal{H}_0$, the joint distribution of $\mathcal{Z}$ is invariant under permutations of the indices, so that for any permutation $\pi$ of $\{1,\dots,n+m\}$, the statistic denoted $\widehat{\MMD}_{k,\pi}^2$ and computed with the sample $(Z_{\pi(1)},\dots,Z_{\pi(n+m)})$ has the same distribution as the original statistic $\widehat{\MMD}_k^2$. One draws $B$ independent permutations $\pi_1,\dots,\pi_B$ uniformly from the set of all permutations of $\{1,\dots,n+m\}$, independently of the data $\mathcal{Z}$, and computes the permuted statistics $\widehat{\MMD}_{k,\pi_b}^2$ for $b=1,\dots,B$. Together with the original value $\widehat{\MMD}_k^2$, this yields $B+1$ exchangeable statistics, and the rejection threshold is defined as the empirical $(1-\alpha)$-quantile of these $B+1$ values \cite{albert2015docteur,schrab2023mmd}. The resulting test has non-asymptotic level $\alpha$ under $\mathcal{H}_0$, a guarantee that follows from the permutation test lemma of Romano and Wolf based on exchangeability arguments \cite{romano2005exact,lehmann2005testing}.

\medskip
Alternative calibration strategies rely on the asymptotic distribution of $\widehat{\MMD}_k^2$ under the null hypothesis, which is given by an infinite weighted sum of independent chi-square random variables \cite{gretton2012kernel}. In practice, this distribution is typically approximated either through spectral methods or by moment-matching with a Gamma distribution \cite{gretton2009fast, gretton2012kernel}. In the remainder of the paper, we consider the permutation procedure, which provides exact finite-sample level control and non-asymptotic theoretical guarantees.

\subsection{Group actions}
\label{sec:group-actions}

We briefly review the basic notions from group theory and group actions used throughout this paper. This discussion is kept concise and emphasizes aspects that are useful for our constructions. For more details on locally compact groups, Haar measures, and integration on groups, see \cite{folland1995,rudin1990}.

\medskip
\textbf{Groups and Haar measure.} A set $G$ equipped with a binary operation $\ast : G \times G \to G$ is a group if $\ast$ is associative (for all $a,b,c$ in $G$, $(a \ast b) \ast c = a \ast (b \ast c)$), there exists an identity element $e \in G$ (for all $a \in G$, $a \ast e = e \ast a = a$) and every element of $G$ has an inverse (for all $a \in G$, there exists $a^{-1} \in G$ such that $a \ast a^{-1} = a^{-1} \ast a = e$). For notational convenience, we write $ab$ instead of $a \ast b$. In the following, $G$ is assumed to be a topological group. That is, a group endowed with a topology for which the group operation $(a,b) \mapsto a \ast b$ and the inversion $a \mapsto a^{-1}$ are continuous. A fundamental result in harmonic analysis states that if $G$ is locally compact and Hausdorff, then there exists a nonzero measure $\lambda_l$ on the Borel $\sigma$-algebra of $G$, called a left Haar measure, which is invariant under left translations. Meaning that, for all $g\in G$ and $A \subseteq G$ a measurable set: $\lambda_l (gA)=\lambda_l (A)$, where $gA = \{ ga \; | \; a \in A \}$. The left Haar measure is unique up to a multiplicative constant. Similarly, there exists a unique right-invariant Haar measure $\lambda_r$. In the sequel, we consider unimodular groups, for which the left and right Haar measures coincide. We then refer to $\lambda$ as the Haar measure, invariant under both left and right translations. Another important result states that $\lambda(G)$ is finite if and only if $G$ is compact. In this case, $\lambda$ can be normalized into a probability measure. When $G$ is non-compact (e.g. $G=\mathbb R^d$), the Haar measure has infinite total mass. Moreover, if $G$ is $\sigma$-compact, the Haar measure is $\sigma$-finite.

\medskip
\textbf{Group actions, orbits, and quotient spaces.}
Let $(\mathcal{X},\mathcal{B}(\mathcal{X}))$ be a measurable space.
A (left) action of $G$ on $\mathcal{X}$ is a mapping
\begin{align*}
\varphi : \; G \times \mathcal{X} \; &\to \; \mathcal{X} \\
(g,x) \; &\mapsto \; \varphi_g(x),
\end{align*}
such that $\varphi_e = id_{\mathcal{X}}$ (the identity function on $\mathcal{X}$) and for all $g,h\in G$, $\varphi_{gh} = \varphi_g \circ \varphi_h$. The orbit of a given $x$ in $\mathcal{X}$ is defined by $[x] :=\{\varphi_g(x) \; | \; g\in G\}$. We denote by $\mathcal{X}/G := \{[x] \mid x \in \mathcal{X}\}$ the set of all orbits, called the quotient space. We also denote by $\Pi$ the canonical projection from $\mathcal{X}$ to $\mathcal{X}/G$, associating to each $x$ in $\mathcal{X}$ its orbit $\Pi (x) := [x]$. By construction, for all $g\in G$, $\Pi\circ \varphi_g=\Pi$. We equip $\mathcal{X}/G$ with the quotient $\sigma$-algebra defined by 
$$\mathcal B (\mathcal{X}/G) :=\{A\subseteq \mathcal{X}/G \; | \; \Pi^{-1}(A)\in\mathcal{B}(\mathcal{X})\},$$
which is the largest $\sigma$-algebra on $\mathcal{X}/G$ that makes $\Pi$ measurable. All probability measures on $\mathcal{X}/G$ considered in this paper are defined on $(\mathcal{X}/G,\mathcal{B}(\mathcal{X}/G))$. Importantly, our theoretical results only rely on this measurable structure and no topological assumptions on $\mathcal{X}/G$ are required. In the following, we assume that the action $\varphi$ is jointly measurable with respect to the product $\sigma$-algebra on $G\times \mathcal{X}$. 

\medskip
\begin{remark}
If $\mathcal{X}$ is a topological space and the action $\varphi$ is continuous, one may endow $\mathcal{X}/G$ with the quotient topology, defined as the finest topology making $\Pi$ continuous.
However, unless additional assumptions (such as properness of the action) are imposed, this topology may be pathological (e.g. non-Hausdorff). Since our analysis is purely measure-theoretic, we do not rely on any topological properties of $\mathcal{X}/G$. For further discussion of transformation groups and quotient spaces (including conditions such as properness ensuring well-behaved quotients), see~\cite{Bredon1972}.
\end{remark}

We present below some examples of the space $\mathcal{X}$, the group $G$ and its associated action.
\begin{itemize}
\item \textbf{Image rotations.} Let $\mathcal{X}$ be a space of images. For example, we may take $\mathcal{X}$ to be a subspace of functions from $\mathbb{R}^2$ to $\mathbb{R}$. Let $G = \mathrm{SO}(2)$ be the special orthogonal group in dimension 2, a compact group, defined by 
\begin{align*}
\mathrm{SO}(2) :=& \left\{ Q \in \mathcal{M}_2(\mathbb{R}) \, | \, Q^\top Q = I_2 \mbox{ and } \det (Q) = 1 \right\} \\
=& \left\{Q_\theta =
\begin{pmatrix}
\cos \theta & -\sin \theta \\
\sin \theta & \cos \theta
\end{pmatrix} 
\in \mathcal{M}_2(\mathbb{R})
\, \middle| \,
\theta \in [0,2\pi)
\right\},
\end{align*}
where $\mathcal{M}_2(\mathbb{R})$ is the set of 2$\times$2 real matrices and $I_2$ is the identity 2$\times$2 matrix. The group action is defined for all $\theta$ in $[0,2\pi)$ and $x$ in $\mathcal{X}$ by 
$$\varphi_\theta (x) : u \mapsto x(Q_\theta^\top u).$$  
The action $\varphi_\theta$ corresponds to rotating the image by an angle $\theta$ around the origin. 

\item \textbf{Time shifts of periodic signals.} Consider $\mathcal{X}$ to be a space of $T$-periodic signals on $\mathbb{R}$, and let $G = \mathbb{R}/T\mathbb{Z}$ act by time shifts. More precisely, the group action is defined for all $[\tau]$ in $G$ and $x$ in $\mathcal{X}$ by
$$\varphi_{[\tau]} (x) : u \mapsto x(u - \tau).$$ 
In this case, $G$ is compact. 

\item \textbf{Time shifts of aperiodic signals.} Let $\mathcal{X}$ be a space of non-periodic signals on $\mathbb{R}$, and let $G = (\mathbb{R},+)$ act on $\mathcal{X}$ by translations. For all $t$ in $G$ and $x$ in $\mathcal{X}$, the action is given by
$$\varphi_t(x) : u \mapsto x(u - t).$$
In this case, $G$ is non-compact but locally compact.\end{itemize}

\section{Invariant two-sample tests under group actions}
\label{sec:invariant-mmd}

The aim of this section is to construct a nonparametric two-sample test for comparing two distributions modulo a group action. For this, we consider two independent i.i.d samples
$\{X_i\}_{i=1}^n \sim P$ and $\{Y_i\}_{i=1}^m \sim Q$, where $P$ and $Q$ are probability measures on a measurable space $(\mathcal{X},\mathcal{A})$, a unimodular group $G$ endowed with its Haar measure $\lambda$, and a group action $\varphi$ of $G$ on $\mathcal{X}$. Under the notations and assumptions of Section \ref{sec:group-actions}, we call a $G$-invariant two-sample test, the following testing problem
\begin{equation}
\label{G-invariant_test}
\mathcal{H}_0 : \Pi_* P = \Pi_* Q \quad \mbox{against} \quad \mathcal{H}_1 : \Pi_* P \neq \Pi_* Q,
\end{equation}
where $\Pi_* P$ (respectively $\Pi_* Q$) denotes the pushforward measure of $P$ (respectively $Q$) by $\Pi$. The idea here is to compare the distributions $P$ and $Q$ while disregarding the variability induced by the group action, which is viewed as a nuisance transformation. While the case of compact groups has been partially studied in the literature \cite{Haasdonk2005,Mroueh2015,Raj2017}, the case of locally compact groups remains unexplored. An objective of this work is to address this extension. In the compact case, the starting point is to consider an averaged kernel 
\begin{equation*}
\ell_\lambda : (x,y) \mapsto \int_G \!\int_G \ell \left(\varphi_g(x),\varphi_{h}(y) \right) \mathrm{d} \lambda(g)\,\mathrm{d} \lambda(h),
\end{equation*}
where $\ell$ is a characteristic kernel on $\mathcal{X}$. This averaging removes the information provided by the group's action and defines a kernel whose values are completely determined by the orbits. For locally compact groups, the kernel $\ell_\lambda$ is not always well-defined, as the Haar measure can be infinite. In the sequel, we assume that $G$ is a locally compact and $\sigma$-compact group and that the two following assumptions hold. The space $\mathcal{X}$ is assumed to be Polish and the group action $\varphi:G\times \mathcal{X}\to \mathcal{X}$ is jointly measurable. We also equip the quotient space $\mathcal{X}/G$ with the quotient $\sigma$-algebra introduced in Section \ref{sec:group-actions}.

\subsection{Weighting and admissible measures}

A crucial step in building the test given by \eqref{G-invariant_test} is to be able to define, for a given probability $P$ on $\mathcal{X}$ and a measure $\nu$ on $G$, an average of the transformations of $P$ under the action of the group $G$. If it is well-defined, this $\nu$-averaged probability is given by
\begin{equation}
\label{Av-Proba}
S_\nu P \; := \; \int_G (\varphi_g)_*P \,\mathrm{d} \nu(g), 
\end{equation}
where $(\varphi_g)_*P$ denotes the pushforward measure of $P$ by $\varphi_g$. We will show later that when $\nu = \lambda$ is the Haar measure, the testing problem \eqref{G-invariant_test} amounts to comparing $S_\nu P$ and $S_\nu Q$. To remain fully general, we define the $\nu$-average of a kernel $\ell$ by
\begin{equation}
\label{Av_Kernel_nu}
\ell_\nu : (x,y) \mapsto \int_G \!\int_G \ell \left(\varphi_g(x),\varphi_{h}(y) \right) \mathrm{d} \nu(g)\,\mathrm{d} \nu(h).
\end{equation}
It is clear that $S_\nu P$ is not always defined in the case of non-compact groups and for a general measure $\nu$. In what follows, we introduce classes of measures on $\mathcal{X}$ and $G$ for which the $\nu$-averaged probability in \eqref{Av-Proba} and the $\nu$-averaged kernel in \eqref{Av_Kernel_nu} are well-defined. For this, we consider a weighting function $\rho$ with a sufficiently fast decay to ensure the integrability with respect to measures on $\mathcal{X}$ and $G$. We assume that $\rho$ is a strictly positive and bounded Borel function and we denote by $\mathcal{M}(\mathcal{X})$ the space of signed measures on $(\mathcal{X},\mathcal B(\mathcal{X}))$ and by $\mathcal{M}_{\sigma}(G)$ the space of $\sigma$-finite measures on $(G,\mathcal B(G))$. We denote by $\mathcal{M}^\rho(\mathcal{X})$ and $\mathcal{M}^\rho(G)$ the following two classes of measures
\begin{align*}
\mathcal{M}^\rho(\mathcal{X})
&:= \left\{ \mu \in \mathcal{M}(\mathcal{X}) : \int_\mathcal{X} \rho(x)\,\mathrm{d}|\mu|(x) < \infty \right\},\\
\mathcal{M}^\rho(G)
&:= \left\{ \nu \in \mathcal{M}_{\sigma}(G) : \sup_{x\in \mathcal{X}} \left\{ \int_G \rho(\varphi_g(x))\,\mathrm{d} \nu (g) \right\}  < \infty \right\}.
\end{align*}

Note that the class $\mathcal{M}^\rho(G)$ includes all probability measures on $G$. More importantly, for a suitable choice of the weighting function $\rho$ and depending on the group action under consideration,
$\mathcal{M}^\rho(G)$ may also contain infinite measures, including the Haar measure $\lambda$ on locally compact and $\sigma$-compact groups.
As stated above, the role of the weighting function $\rho$ is precisely to control integrability along group orbits, thereby extending the class of admissible measures on $G$ beyond finite measures. Consider now $k$, a continuous, bounded and positive definite kernel on $\mathcal{X}$ and the weighted kernel 
$$k^{\rho} : (x,y) \mapsto \rho(x) k (x,y) \rho(y).$$
Since $k$ is measurable and positive definite, and $\rho$ is measurable, it follows that the kernel $k^{\rho}$ is measurable and positive definite on $\mathcal{X} \times \mathcal{X}$. With these definitions in place, we now provide a sufficient condition ensuring that the $\nu$-averaged probability $S_\nu P$ is well-defined. We then establish a result showing that the distance induced by the weighted kernel $k^\rho$ distinguishes $\nu$-averaged probabilities. 
\begin{prop}
\label{prop:S_nu}
Let $\nu\in \mathcal{M}^\rho(G)$. For all probability measures $P \in \mathcal P(\mathcal{X})$, the averaged measure
$$S_\nu P := \int_G (\varphi_g)_* P\, \mathrm{d}\nu(g)$$
is well-defined and belongs to $\mathcal{M}^\rho(\mathcal{X})$.
\end{prop}
\begin{prop}
\label{prop:carac}
Assume that $k$ is characteristic on the space of finite signed measures $\mathcal{M}_f(\mathcal{X})$. Then, $k^\rho$ is characteristic on $\mathcal{M}^\rho(\mathcal{X})$.
\end{prop}
This means that the MMD based on $k^\rho$ is able to distinguish averaged distributions, which is a crucial ingredient for consistency of the invariant test. To estimate $\MMD_{k^\rho} (S_\nu P, S_\nu Q)$ from i.i.d. samples drawn from $P$ and $Q$, we first express it as an MMD between $P$ and $Q$ with respect to another kernel. It turns out that this kernel is precisely the $\nu$-average of the kernel~$k^\rho$.

\subsection{Averaged and invariant kernels}

Let us introduce the averaged kernel that will be used to perform the $G$-invariant two-sample test defined in \eqref{G-invariant_test}. The main idea is to build a kernel that compares observations only through their behavior along group orbits. This is achieved by averaging the weighted kernel $k^\rho$ over the action of the group. Let $\nu \in \mathcal{M}^\rho(G)$, we define the $\nu$-averaged kernel $k^\rho_\nu$ as 
$$k^\rho_\nu : (x,y)
\mapsto \int_G \!\int_G k^\rho \left(\varphi_g(x),\varphi_{h}(y)\right) \mathrm{d} \nu(g)\,\mathrm{d} \nu(h).$$
The integrability conditions encoded in the definition of $\mathcal{M}^\rho(G)$ ensure that the kernel $k^\rho_\nu$ is well-defined. Indeed, for $x,y$ in $\mathcal{X}$
\begin{align}
k^\rho_\nu (x,y) &= \int_G \!\int_G k^\rho \left(\varphi_g(x),\varphi_{h}(y)\right) \mathrm{d} \nu(g)\,\mathrm{d} \nu(h) \nonumber \\ 
&= \int_G \!\int_G \rho\left(\varphi_g(x) \right) k \left(\varphi_g(x),\varphi_{h}(y)\right) \rho\left(\varphi_{h}(y) \right) \mathrm{d} \nu(g)\,\mathrm{d} \nu(h) \nonumber\\
&\leq  \Vert k \Vert_{\infty} \left[ \int_G \rho\left(\varphi_g(x) \right) \mathrm{d} \nu(g)  \right] \left[ \int_G \rho\left(\varphi_{h} (y)\right) \mathrm{d} \nu(h) \right]  \nonumber \\
&< +\infty \label{Av-ker-finite}
\end{align}
where $\Vert k \Vert_{\infty} = \sup_{x,y \in \mathcal{X}} |k(x,y)|$. The kernel $k^\rho_\nu$ compares two elements $x$ and $y$ by averaging similarities between pairs of points taken from their respective orbits. If $k$ is symmetric positive definite, then $k^\rho_\nu$ is as well. This result is shown in the next proposition, which also provides a useful interpretation of the feature representation of $k^\rho_\nu$.
\begin{prop}
\label{prop:Av-kernel}
Assume that $k$ and $\rho$ are continuous and bounded. Let $\nu$ in $\mathcal{M}^\rho(G)$ and $\Phi^\rho : x \mapsto k^\rho(x,\cdot)$ be the canonical feature map of $k^\rho$. Then, for all $x,y$ in $\mathcal{X}$,
$$k^\rho_\nu(x,y)
= \left\langle \Phi^\rho_\nu(x), \Phi^\rho_\nu(y)\right\rangle_{\mathcal H^\rho_{k}},$$
where $\displaystyle \Phi^\rho_\nu : x \mapsto \int_G \Phi^\rho \left( \varphi_g (x) \right) \mathrm{d} \nu(g)$. In other words, $\Phi^\rho_\nu$ is the canonical feature map of $k^\rho_\nu$. In particular, $k^\rho_\nu$ is a positive definite kernel on $\mathcal{X}$.
\end{prop}

\medskip
Proposition~\ref{prop:Av-kernel} shows that the feature map associated with $k^\rho_\nu$ is obtained by averaging the one of $k^\rho$ along group orbits. The natural question that arises is the link between the mean embeddings associated with the kernels $k^\rho$ and $k_\nu^\rho$. This result is given in the next proposition. 
\begin{prop}
\label{prop:Av-mean-emb}
Assume that $k$ and $\rho$ are continuous and bounded. Let $\nu$ be in $\mathcal{M}^\rho(G)$ and $P$ in $\mathcal P(\mathcal{X})$. Denote by $m^\rho_\nu (P)$ the mean embedding of $P$ with respect to $k^\rho_\nu$, and by $m^\rho (S_\nu P)$ the mean embedding of $S_\nu P$ with respect to $k^\rho$. Then, 
$$m^\rho_\nu (P) = m^\rho (S_\nu P).$$
\end{prop}
Proposition \ref{prop:Av-mean-emb} shows that embedding $P$ with $k^\rho_\nu$ is equivalent to embedding the averaged measure $S_\nu P$ with $k^\rho$. This can also be interpreted in light of Proposition \ref{prop:Av-kernel}. The philosophy behind $m^\rho_\nu (P)$ is to embed first and then average, whereas that of $m^\rho (S_\nu P)$ is to average first and then embed. These are two sides of the same coin. The next corollary is an immediate consequence.
\begin{cor}
\label{cor:MMD-k-knu}
Assume that $k$ and $\rho$ are continuous and bounded. Let $\nu$ in $\mathcal{M}^\rho(G)$ and $P,Q$ in $\mathcal P(\mathcal{X})$. Then,  
$$\MMD_{k^\rho_\nu} (P,Q)  = \MMD_{k^\rho} (S_\nu P, S_\nu Q).$$
\end{cor}
Recall that $k^\rho$ is characteristic on $\mathcal{M}^\rho(\mathcal{X})$ and, by Proposition~\ref{prop:S_nu}, that $S_\nu P \in \mathcal{M}^\rho(\mathcal{X})$ for every probability distribution $P$. It follows that $\MMD_{k^\rho} (S_\nu P, S_\nu Q) = 0$ if and only if $S_\nu P = S_\nu Q$.

\medskip
Let us now focus on the case $\nu = \lambda$, where $\lambda$ denotes the Haar measure. When $G$ is compact, this measure is finite and can therefore be normalized to a probability measure, namely the uniform distribution on $G$. In this setting, without further assumptions on $\lambda$, both the $\lambda$-averaged probability $S_\lambda P$ and the $\lambda$-averaged kernel $k^\rho_\lambda$ in \eqref{Av-Proba} and \eqref{Av_Kernel_nu} are well-defined. When $G$ is non-compact, the Haar measure has infinite total mass. However, for a suitable choice of the weighting function $\rho$, it belongs to $\mathcal{M}^\rho(G)$. Then, according to Proposition \ref{prop:S_nu} and Equation~\eqref{Av-ker-finite}, the averaged probabilities and the averaged kernel are well-defined, despite the infiniteness of $\lambda$.

\begin{theorem}
\label{prop:Haar-quotient}
Assume that $k$ and $\rho$ are continuous and bounded, and that the Haar measure $\lambda$ belongs to $\mathcal{M}^\rho(G)$. Let $P,Q$ in $\mathcal{P} (\mathcal{X})$. Then, $S_\lambda P = S_\lambda Q$ if and only if $\Pi_* P = \Pi_* Q$. 
\end{theorem}

The distribution $S_\lambda P$ is obtained by averaging $P$ uniformly along each group orbit. Because the Haar measure is invariant, this averaging treats all points in the orbit in the same way. As a result, any information about how the mass is distributed inside an orbit disappears, and only the total mass assigned to each orbit remains. Since $\Pi_* P$ exactly represents this mass on the orbit space, $S_\lambda P$ is completely determined by $\Pi_* P$ and conversely. If $\nu$ is not the Haar measure, the averaging is no longer uniform along the orbits. In that case, the result may still depend on how the mass is arranged within each orbit, and this correspondence with the quotient distribution no longer holds. 
\begin{remark}
\label{Rq:rho}
The condition $\lambda$ in $\mathcal{M}^\rho(G)$ depends on both the action and the choice of $\rho$. In general, there is no universal choice of $\rho$ ensuring it. Consider the case where $\varphi$ is the trivial action, namely $\varphi : (g,x) \mapsto x$. If $G$ is non-compact, then any weighting function $\rho > 0$ satisfies
$$\int_G \rho \left( \varphi_g (x) \right) \mathrm{d}\lambda(g) = \int_G \rho(x) \mathrm{d} \lambda(g) = \rho(x)\lambda(G) = +\infty.$$
This means that $\lambda$ does not belong to $\mathcal{M}^\rho(G)$. At the opposite extreme, when $G$ is compact, the Haar measure is finite. One may simply take $\rho\equiv 1$, in which case $k^\rho_\lambda$ yields the classical Haar-averaged invariant kernel.    
\end{remark}
It is also important to emphasize that since $k^\rho_\lambda$ is constant on orbits, it induces a measurable positive kernel $\widetilde{k}^\rho_\lambda$ on the quotient $\mathcal{X}/G$, defined by
\begin{equation}
\label{kernel-orbits}
\widetilde{k}^\rho_\lambda : \left( [x],[y] \right) \mapsto k^\rho_\lambda (x,y).
\end{equation}
The kernel $\tilde{k}_\lambda$ is characteristic on the set of pushforward measures under $\Pi$, namely
$$\Pi_* \mathcal{P}(\mathcal{X}) := \left\{ \Pi_*P \; | \; P\in\mathcal \mathcal{P}(\mathcal{X}) \right\}.$$

\medskip
Thanks to Proposition~\ref{prop:Haar-quotient} and Corollary \ref{cor:MMD-k-knu}, by using the invariant kernel $k^\rho_\lambda$, one can perform the MMD two-sample test on the quotient space without ever having to construct it explicitly. In other words, the testing problem \eqref{G-invariant_test} boils down to testing
\begin{equation}
\label{G-invariant_test_MMD}
\mathcal{H}_0 : \MMD_{k^\rho_\lambda} (P,Q) = 0 \quad \mbox{against} \quad \mathcal{H}_1 : \MMD_{k^\rho_\lambda} (P,Q) > 0.
\end{equation}
In the next section, we discuss the practical aspects of its implementation.

\subsection{Invariant MMD test in practice}
\label{subsec:application-mmd}

Assuming that the kernel $k^\rho_\lambda$ is known, a two-sample test can be performed based on \eqref{G-invariant_test_MMD}. To do so, the MMD is estimated using the unbiased U-statistic introduced in Section \ref{sec:kernel_review}. The rejection threshold is then obtained via the permutation procedure. In our setting, the theoretical validity of this procedure still holds. Indeed, the kernel $k^\rho_\lambda$ depends only on the orbits through the kernel $\widetilde{k}^\rho_\lambda$ introduced in \eqref{kernel-orbits}. Moreover, under $\mathcal{H}_0$, the pooled sample of orbits $\left\{\Pi(X_1),\dots,\Pi(X_n),\Pi(Y_1),\dots,\Pi(Y_m) \right\}$ is exchangeable. Consequently, the distribution of the MMD estimator is unchanged under permutations. Therefore, the permutation procedure yields a valid level-$\alpha$ test.

\medskip
Practically speaking, computing $k^\rho_\lambda$ requires evaluating an integral that is generally intractable in closed form. To address this issue, we introduce a procedure to approximate it for each pair $(x,y) \in \mathcal{X}^2$. Assume that the Haar measure $\lambda$ belongs to $\mathcal{M}^\rho(G)$ for some weighting function~$\rho$. We define the orbit-averaged weight function by
\begin{equation*}
\label{eq:app-w-def}
w : x \mapsto \int_G \rho \left( \varphi_g (x) \right) \mathrm{d}\lambda(g).
\end{equation*}
For all $x$ in $\mathcal{X}$, define the measure $\nu_x$ on $G$ by
\begin{equation}
\label{eq:app-nux-def}
\nu_x(dg) := \frac{\rho \left( \varphi_g (x) \right)}{w(x)}\,\mathrm{d} \lambda(g).
\end{equation}
Then, $\nu_x$ is a probability measure on $G$. In addition, $k^\rho_\lambda$ can be written as 
\begin{equation}
\label{eq:app-klambda-exp}
k^\rho_\lambda : (x,y)
\mapsto w(x) \, w(y) \int_G \!\int_G k \left(\varphi_g (x),\varphi_h (y) \right) \mathrm{d} \nu_x (g) \,  \mathrm{d} \nu_y (h).
\end{equation}
The representation \eqref{eq:app-klambda-exp} plays a key role in approximating the kernel $k^\rho_\lambda$ for pairs of points $(x,y)$ in $\mathcal{X}^2$. Indeed, it enables a Monte-Carlo approximation. We assume that $\rho$ is chosen so that the normalizing constant $w(x)$ is analytically computable or numerically approximable. Note that when $G$ is compact, Remark~\ref{Rq:rho} implies that $\rho \equiv 1$, $w \equiv \lambda(G)$ and that for all $x$ in $\mathcal{X}$, $\nu_x = \lambda/\lambda(G)$. Let $\{g_s\}_{s=1}^S$ and $\{h_s\}_{s=1}^S$ be i.i.d. samples from $\nu_x$ and $\nu_y$, respectively. Then, $k^\rho_\lambda (x,y)$ can be approximated by 
$$\overline{k}^\rho_\lambda (x,y)
:= \frac{\widetilde{w} (x) \widetilde{w} (y)}{S^2} \sum_{r=1}^S \sum_{s=1}^S k \left( \varphi_{g_r} (x) ,\varphi_{h_s} (y) \right),$$
where $\widetilde{w} (x)$ and $\widetilde{w} (y)$ are approximations of $w(x)$ and $w(y)$, respectively. Once the kernel $k^\rho_\lambda$ is approximated by $\overline{k}^\rho_\lambda$, the permutation-based test described in Section~\ref{sec:kernel_review} can be carried out using samples from $P$ and $Q$. Obviously, increasing $S$ reduces the approximation error, but results in higher computational cost and memory usage. We now present two illustrative examples concerning the choice of $\rho$ and the approximation of $k^\rho_\lambda$. Recall that all the theoretical results presented in this paper are stated in a very general setting. The only assumption on $\mathcal{X}$ is that it is a Polish space. In particular, it is not required to be locally compact. This generality allows us to consider applications in functional data analysis, where observations are curves or signals. The Hilbert space $\mathcal{X} = L^2 (I)$ of real-valued square-integrable functions on an interval $I \subseteq \mathbb{R}$, provides a natural and widely used functional framework. A commonly used kernel on $L^2 (I)$ is the Gaussian one, defined by
\begin{equation}
\label{eq:gaussian-L2}
k: (x,y) \mapsto \exp \left(-\frac{\|x-y\|_2^2}{2\sigma^2} \right),
\end{equation}
where $\| \cdot \|_2$ denotes the usual $L^2$-norm and $\sigma>0$ is a bandwidth parameter. It has recently been shown in \cite{wynne2022kernel} that the Gaussian kernel is characteristic on the class of finite measures. Therefore, it constitutes a sound and flexible choice for the base kernel.

\medskip
\textbf{Time shifts of periodic signals.} Consider the space of real-valued $1$-periodic functions on $\mathbb{R}$ whose restriction to $[0,1]$ belongs to $L^2([0,1])$, modulo equality almost everywhere. This functional space is canonically identified with $\mathcal{X} = L^2(\mathbb{S}^1)$, where $\mathbb{S}^1 = \mathbb{R}/\mathbb{Z}$ denotes the circle. Let $G = (\mathbb{S}^1, +)$ and let $\varphi$ be the circular shift action 
\begin{align*}
\varphi : \; G \times \mathcal{X} \; &\to \; \mathcal{X} \\
(\tau,x) \; &\mapsto \; \varphi_\tau (x) = x(\cdot - \tau).
\end{align*}
As mentioned earlier, since $G$ is compact, we have $\rho \equiv 1$. Furthermore, $\lambda$ can be identified with the Lebesgue measure on $[0,1)$, so that $\lambda(G) = 1$ and $\nu_x$ is the uniform distribution on $\mathbb{S}^1$. Then, the Haar-averaged invariant kernel is defined by
$$
k^\rho_\lambda : (x,y)
\mapsto \int_{\mathbb{S}^1} \!\int_{\mathbb{S}^1} k \left(\varphi_\tau(x),\varphi_{\iota}(y)\right) \mathrm{d} \lambda(\tau)\,\mathrm{d} \lambda(\iota),
$$
where $k$ is the Gaussian kernel introduced in \eqref{eq:gaussian-L2}. Let $\{\tau_s\}_{s=1}^S$ and $\{\iota_s\}_{s=1}^S$ be i.i.d. samples drawn uniformly from $\mathbb{S}^1$. Then, $k^\rho_\lambda (x,y)$ is approximated by 
\begin{equation}
\label{Eq:Approx_Aver_Ker}
\overline{k}^\rho_\lambda (x,y)
= \frac{1}{S^2} \sum_{r=1}^S \sum_{s=1}^S  \exp \left(-\frac{\|x(\cdot - \tau_r) - y(\cdot - \iota_s) \|_2^2}{2\sigma^2} \right).
\end{equation}
In practice, the $L^2$-norms are replaced by discretized approximations.

\medskip
\textbf{Time shifts of aperiodic signals.} Consider now the non-periodic setting where $$\mathcal{X} = \left\{ f \in L^2(\mathbb R) \; | \; \Vert f \Vert^2_2 \leq R \right\},$$ 
where $R > 0$. Let $G = (\mathbb R,+)$ and let $\varphi : (\tau, x) \mapsto x (\cdot - \tau)$ be the time shift action. Here, $G$ is non-compact and its Haar measure $\lambda$ is the Lebesgue measure. For $c > 0$, define the weighting function $\rho_c$~by
\begin{equation}
\label{eq:rho-gauss-window}
\rho_c : x \mapsto
\int_{\mathbb R} x^2(u) \exp \left(-\frac{u^2}{2 c^2} \right) \mathrm{d}u.
\end{equation}
This functional applies a Gaussian window centered at zero, giving more weight to values of the signal near the reference time origin. Let $c > 0$ and $x$ in $\mathcal{X}$. Then, using Tonelli's theorem, we have
\begin{align*}
\int_{\mathbb R} \rho_c \left( \varphi_\tau (x) \right) \mathrm{d}\lambda(\tau) &= \int_{\mathbb R} \int_{\mathbb R} x^2(u - \tau) \exp \left(-\frac{u^2}{2 c^2} \right) \mathrm{d}u \, \mathrm{d}\tau \\
&= \int_{\mathbb R} \left[ \int_{\mathbb R} x^2(u - \tau) \mathrm{d}\tau \right] \exp \left(-\frac{u^2}{2 c^2} \right) 
\mathrm{d}u \\
&= c \sqrt{2\pi} \Vert x \Vert^2_2.  
\end{align*}
Then,
$$\sup_{x\in \mathcal{X}} \left\{ \int_\mathbb{R} \rho_c \left( \varphi_\tau(x) \right) \,\mathrm{d} \lambda (\tau) \right\} \leq c \sqrt{2\pi} R.$$
In other words, $\lambda$ belongs to $\mathcal{M}^{\rho_c}(G)$. Assume now that the function $x$ is discretized over an interval $[-T/2,T/2]$ using a uniform grid with $p$ points $x(t_1), \ldots, x(t_p)$. The discretized version of $\rho_c \left( \varphi_\tau(x) \right)$ is given by
$$\widetilde{\rho}_c \left( \varphi_\tau(x) \right) = \frac{T}{p-1} \sum_{k=1}^p x^2(t_k) \exp \left(-\frac{t_k^2}{2 c^2} \right).$$
This shows that $\nu_x$ can be approximated by a mixture Gaussian distribution. In addition, $w(x)$ can be approximated by 
$$\widetilde{w} (x) = \sqrt{2\pi} c \frac{T}{p-1} \sum_{k=1}^p x^2(t_k).$$

\section{Simulations}
\label{sec:simulations}

We present a numerical study illustrating the implementation of the invariant MMD test developed in this paper. To this end, we consider $\mathcal{X} = L^2 (I)$, where $I$ is a real interval. The group action corresponds to translations along the horizontal axis. Our goal is twofold. On the one hand, we empirically verify level control when the two distributions differ only through the group action. On the other hand, we assess the power when the difference cannot be explained by a translation. We consider two different settings, a periodic one in Section \ref{sec:simu_periodic} and an aperiodic one in Section \ref{sec:simu_aperiodic}. In both settings, we compare the performance of the invariant test, the align-then-test procedure, and the base-kernel test.

\subsection{General setup}

For all kernel-based tests, the test statistic is the unbiased U-statistic MMD estimator recalled in Section~\ref{sec:kernel_review}. In each setting, we simulate i.i.d. samples $\mathbb{X}_n = (X_1,\dots,X_n)$ and $\mathbb{Y}_m = (Y_1,\dots,Y_m)$ drawn from $P$ and $Q$, respectively, with $n=m=20$. Each experiment is repeated $N_{\mathrm{rep}}=300$ times in order to estimate empirical rejection rates. The realizations are discretized on a uniform grid with $p=128$ points over an interval $[a,b]$. In addition, the squared $L^2$-norm of a curve $x$ is approximated as
$$\| x \|_2^2 \approx \frac{b-a}{p-1} \sum_{k=1}^p x^2(t_k).$$
On a discrete grid, the action $\varphi_t z = z(\cdot - t)$ is implemented by evaluating the signal at shifted locations $t_k-t$ using linear interpolation (with periodic wrapping in Experiment~1 and zero padding in Experiment~2).
In the main simulations, the kernel $k^\rho_\lambda$ is approximated by its Monte Carlo version $\overline{k}^\rho_\lambda$ according to Equation \eqref{Eq:Approx_Aver_Ker}, with sample sizes $S=16$. This provides a good trade-off between computational cost and statistical performance. Additional numerical simulations are conducted with $S = 32$. The base kernel $k$ is the Gaussian kernel defined in Equation \eqref{eq:gaussian-L2}. The bandwidth $\sigma$ is chosen at each repetition using the median heuristic applied to the pairwise $L^2$-distances of the pooled sample $\mathbb{Z}_n=(X_1,\dots,X_n,Y_1,\dots,Y_m)$. The rejection threshold is calibrated via the permutation procedure described in Section \ref{sec:kernel_review}. To do so, we generate uniformly $B =200$ random permutations $\{\pi_b\}_{b=1}^{B}$ of the labels, compute the permuted statistics $\widehat{\mathrm{MMD}}^{2}_{k,\pi_b}$ and estimate the $p$-value as
$$p_{\mathrm{val}} = \frac{1}{B+1} \left( 1 + \sum_{b=1}^{B} \mathbf{1}_{ \left\{ \widehat{\mathrm{MMD}}^{2}_{k,\pi_b} \geq \widehat{\mathrm{MMD}}^{2}_{k} \right\}} \right),$$
where $\widehat{\mathrm{MMD}}^{2}_{k}$ is the U-statistic computed with the original sample $\mathbb{Z}_n$. Then, we reject the null hypothesis whenever $p_{\mathrm{val}} \le \alpha$ with $\alpha=0.05$.

\medskip
\textbf{Alignment baseline.} Besides comparing the kernels $k$ and $k^\rho_\lambda$, we also consider the classical baseline "align-then-test". This approach attempts to remove spurious translations by aligning the curves before performing the test. To do so, one defines a reference signal for each sample, after which the observations are shifted and scaled to best match it. We denote by $X_{\mathrm{ref}}$ and $Y_{\mathrm{ref}}$ the reference signals associated with $\mathbb{X}_n$ and $\mathbb{Y}_m$, respectively. In the remainder, these reference signals are considered to be the empirical medoids. For the sample $\mathbb{X}_n$, the reference is defined by 
$$X_{\mathrm{ref}} \in \underset{X_i \in \mathbb{X}_n}{\arg\min} \sum_{X_j \in \mathbb{X}_n} \|X_i - X_j\|_2.$$
Then, each realization $X_i$ is shifted by $\hat{t}_i$ given by 
$$\hat{t}_i \in \underset{t}{\arg\min} \inf_{a} \|X_{\mathrm{ref}} - a X_i (\cdot - t) \|_2.$$
A similar definition holds for $Y_{\mathrm{ref}}$ and the same alignment procedure is applied. Then, the MMD two-sample test is performed using the aligned samples and the kernel $k$. Such alignment pipelines are standard in functional data analysis and signal processing (see, e.g., \cite{ramsay2005functional}). Unlike the invariant-kernel approach, this method requires estimating nuisance parameters and selecting a reference signal. This may affect its stability when alignment is weakly identifiable.

\subsection{Time shifts of periodic signals}
\label{sec:simu_periodic}

We consider a periodic setting where the observations are $2\pi$-periodic, up to multiplicative noise. The actions correspond to circular shifts on $[0,2\pi)$. As in the periodic example of Section \ref{subsec:application-mmd}, the acting group $G = \mathbb{R}/2\pi \mathbb{Z}$ is compact and we take $\rho \equiv 1$. 

\medskip
\textbf{Data generation.} Let $X$ and $Y$ be the random processes defined by 
$$X = \gamma h_1( \cdot -\theta_1)\varepsilon \quad \text{and} \quad Y =  \gamma h_2(\cdot -\theta_2)\varepsilon,$$
where the amplitude $\gamma$ is log-normally distributed with log-mean 0 and log-standard deviation $\sigma_\gamma$, $\theta_1$ and $\theta_2$ are random phases (taken modulo $2\pi$) and $\varepsilon$ is Gaussian white noise with mean~1 and standard deviation 0.8. The random variables $\gamma$, $\theta_1$, $\theta_2$ and the random process $\varepsilon$ are independent. The distributions and the functions $h_1$, $h_2$ vary across simulations. They are specified in the corresponding scenarios below. The sample $\mathbb{X}_n$ (respectively $\mathbb{Y}_n$) is generated according to the distribution of $X$ (respectively $Y$).

\medskip
\textbf{Scenario $\mathcal{H}_0$.} In this scenario, we consider sinusoidal functions $h_1 = h_2 = \sin$, 
$$\theta_1 \equiv \widetilde{\theta}_1 \pmod{2\pi} \quad \mbox{ and } \quad \theta_2 \equiv \widetilde{\theta}_2 \pmod{2\pi},$$
where $\widetilde{\theta}_1$ (respectively $\widetilde{\theta}_2$) is Gaussian with mean $\delta/2$ (respectively $- \delta/2$) and standard deviation $\sigma_{\theta_1} = 0.8$ (respectively $\sigma_{\theta_2} = 0.8$). All values of $\delta$ ranging from 0 to 1 with a step of 0.2 are considered in the simulations. With this generating mechanism, the distributions $P$ and $Q$ differ, but coincide modulo translations. In other words, $\Pi_*P = \Pi_*Q$.  

\medskip
\textbf{Scenario $\mathcal{H}_1$.} Let $h_1 = \sin$, $h_2 : t \mapsto h_1(t) + \delta \sin(2t+0.3)$,
$$\theta_1 \equiv \widetilde{\theta}_1 \pmod{2\pi} \quad \mbox{ and } \quad \theta_2 \equiv \widetilde{\theta}_2 \pmod{2\pi},$$
where $\widetilde{\theta}_1$, $\widetilde{\theta}_2$ are Gaussian with mean 0 and standard deviation $\sigma_{\theta} = 0.8$. In this setting, the distributions differ even after alignment. This means that, $\Pi_*P \neq \Pi_*Q$. 

\begin{figure}[h]
    \centering
    \includegraphics[width=0.49\textwidth]{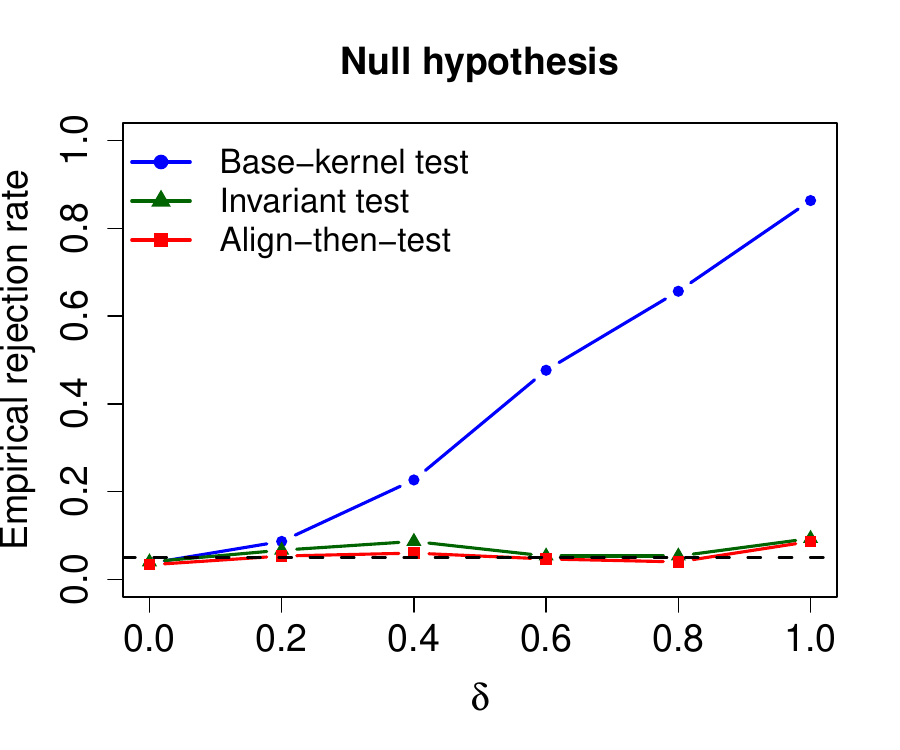}
    \hfill
    \includegraphics[width=0.49\textwidth]{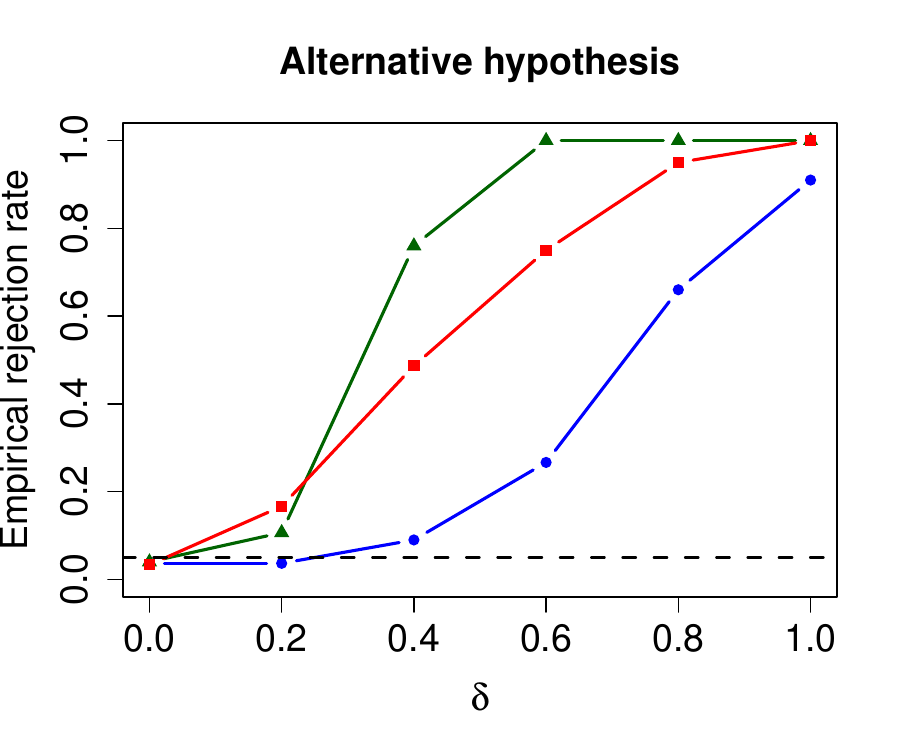}
    \caption{\textbf{Periodic case}. Empirical rejection rates of the MMD tests based on $k$, $k^\rho_\lambda$ and the align-then-test baseline with respect to $\delta$. The approximation budget of $k^\rho_\lambda$ is $S=16$. \textbf{Left.} $P \neq Q$ and $\Pi_*P = \Pi_*Q$. \textbf{Right.} $\Pi_*P\neq\Pi_*Q$.}
    \label{fig:exp1-periodic}
\end{figure}

\begin{figure}[h]
    \centering
    \includegraphics[width=0.48\textwidth]{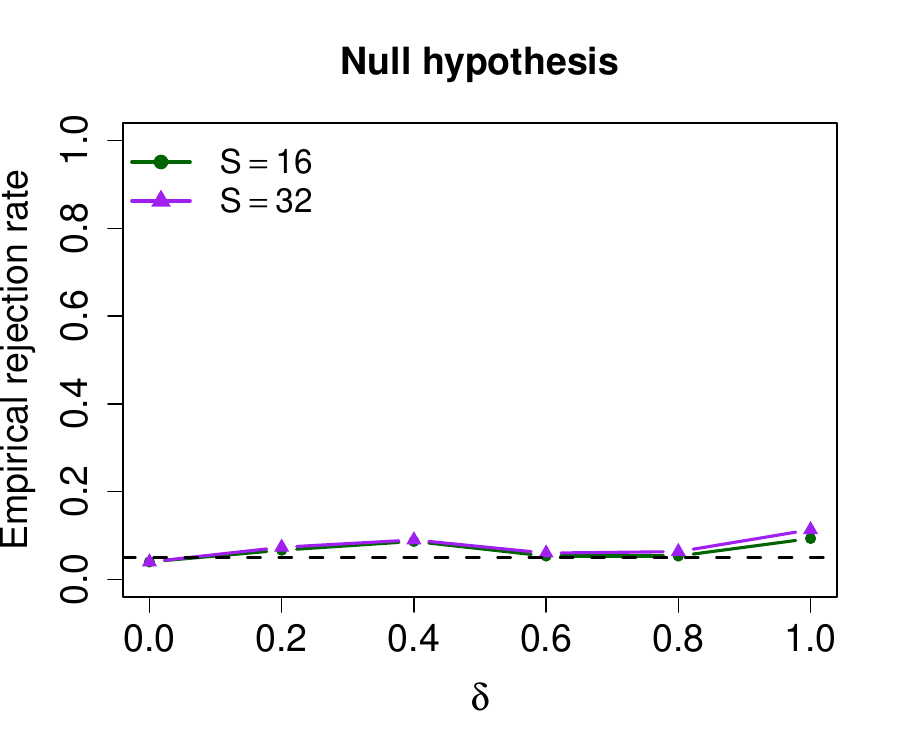}
    \hfill
    \includegraphics[width=0.48\textwidth]{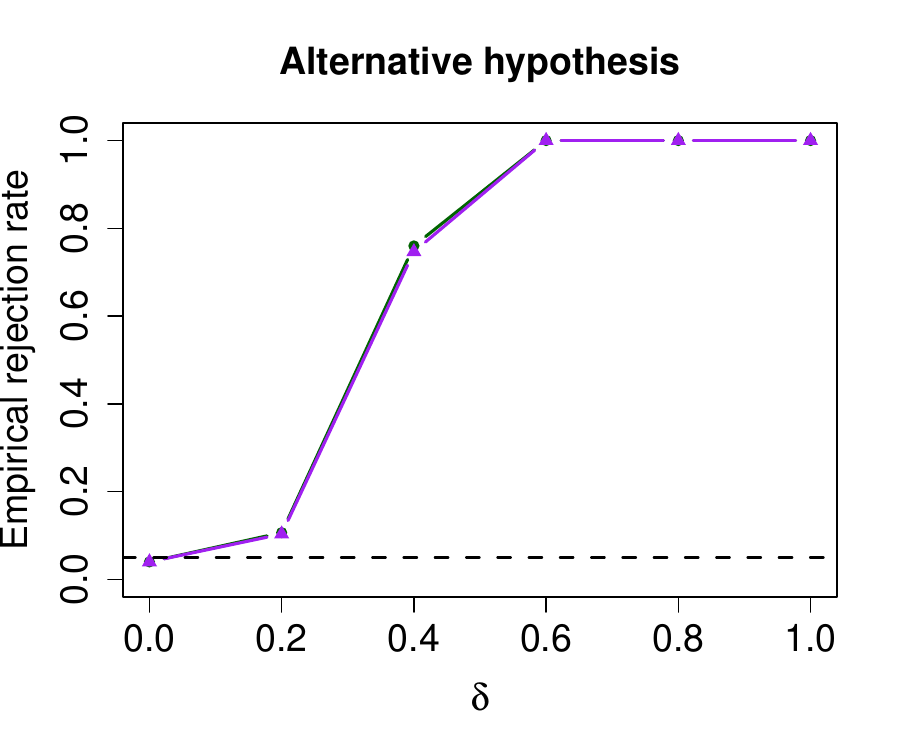}
    \caption{\textbf{Periodic case}. Empirical rejection rates of the MMD tests based on $k^\rho_\lambda$with two different approximation budgets. \textbf{Left.} $P \neq Q$ and $\Pi_*P = \Pi_*Q$. \textbf{Right.} $\Pi_*P\neq\Pi_*Q$.}
    \label{fig:exp1-periodic-S32}
\end{figure}

\medskip
Figure \ref{fig:exp1-periodic} represents the empirical rejection rates, with respect to $\delta$ under the $\mathcal{H}_0$ and $\mathcal{H}_1$ scenarios. Under $\mathcal{H}_0$, the invariant and alignment tests behave well, with low rejection rates. Both procedures are correctly calibrated, their empirical rejection rates remain close to the nominal level $\alpha$ for all values of $\delta$. The results obtained with the kernel $k$ are consistent with its theoretical construction. Since the original distributions differ, the test based on $k$ detects differences and its power increases with $\delta$. Under the alternative, the three tests have increasing power as $\delta$ increases. Moreover, both the invariant test and the align-then-test procedures achieve substantially higher power than the $k$-based test. Therefore, removing the translation nuisance enhances the ability to capture intrinsic shape differences between the two groups of signals. In addition, except for $\delta = 0.2$, the invariant test has significantly higher power than the align-then-test baseline. We also assess the sensitivity of the results to the approximation of $k^\rho_\lambda$. For this, we conduct the same experiments with a larger Monte Carlo budget, namely $S=32$. Figure \ref{fig:exp1-periodic-S32} presents the corresponding empirical rejection rates. The results remain quasi-unchanged. This indicates that the behavior of the invariant test is stable with respect to the approximation budget.

\subsection{Time shifts of aperiodic signals}
\label{sec:simu_aperiodic}

We now consider an aperiodic setting in which the signals are observed on a finite window, namely $[-5,5]$. In this case, the acting group $(\mathbb{R},+)$ is non-compact and a weighting function is required. We choose a weighting function of the form $\rho_c$ defined in Equation \eqref{eq:rho-gauss-window}, where the selection rule for $c$ is introduced below. 

\medskip
\textbf{Choice of the parameter $c$.} We propose the following heuristic, which appears to be robust in practice. Let $\mathbf{t}=(t_1,\dots,t_p)$ be the temporal grid. For a given signal $x$, we define the weighted temporal mean $\mu_x$ and variance $s_x^2$ as
$$\mu(x) = \frac{\sum_{k=1}^p t_k x^2(t_k)}{\sum_{k=1}^p x^2(t_k)} \quad \mbox{and} \quad
s(x)^2 = \frac{\sum_{k=1}^p (t_k - \mu_x)^2 x^2(t_k)}{\sum_{k=1}^p x^2(t_k)}.$$
We then set $c$ to 
$$c := \mathrm{median} \left\{s_i \, | \, 1 \leq i \leq n + m ,\mbox{ and }  s_i > 0 \right\},$$
where $s_i = s(X_i)$ for $1 \leq i \leq n$ and $s_i = s(Y_i)$ for $n + 1 \leq i \leq n + m$. This rule can be interpreted as follows. The quantity $s(x)$ measures the temporal dispersion of the signal energy around its barycenter. The median is then computed over the pooled sample. This provides a robust scale shared by the two groups. As a result, the Gaussian weighting is neither overly localized nor overly diffuse.

\medskip
\textbf{Data generation.} Let $X$ and $Y$ be the random processes defined by
$$X = \gamma h_1 (\cdot - \theta_1) \varepsilon \quad \text{and} \quad Y = \gamma h_2 (\cdot - \theta_2) \varepsilon,$$
where the amplitude $\gamma$ is log-normally distributed with log-mean 0 and log-standard deviation $\sigma_\gamma$, $\theta_1$ and $\theta_2$ are random translations and $\varepsilon$ is Gaussian white noise with mean~1 and standard deviation 0.8. The random variables $\gamma$, $\theta_1$, $\theta_2$ and the random process $\varepsilon$ are independent. As in Section \ref{sec:simu_periodic}, the functions and parameters vary across simulations. The samples $\mathbb{X}_n$ and $\mathbb{Y}_n$ are generated according to the distribution of $X$ and $Y$, respectively.

\medskip
\textbf{Scenario $\mathcal{H}_0$.} In this case, $h_1 = h_2 : t \mapsto \exp (- 2 t^2)$ and $\theta_1$ (respectively $\theta_2$) is Gaussian with mean $\delta/2$ (respectively $-\delta/2$) and standard deviation $\sigma_{\theta_1} = 0.8$ (respectively $\sigma_{\theta_2} = 0.8$). The simulations consider values of $\delta$ ranging from 0 to 1 with increments of 0.2. Similarly to the $\mathcal{H}_0$ case in Section \ref{sec:simu_periodic}, the distributions $P$ and $Q$ are different, but induce the same distribution in the quotient space. 

\medskip
\textbf{Scenario $\mathcal{H}_1$.} Let $h_1 : t \mapsto \exp (- 2 t^2)$ and $h_2 = h_1 + \delta/4 \times p$, where $p$ is the two-bump function defined by 
$$p : t \mapsto \exp \left( -2 (t-1)^2  \right) + 0.4 \exp \left(-\frac{(t+1)^2}{2} \right).$$
The random variables $\theta_1$, $\theta_2$ are Gaussian with mean 0 and standard deviation $\sigma_{\theta} = 0.8$.

\begin{figure}[H]
    \centering
    \includegraphics[width=0.49\textwidth]{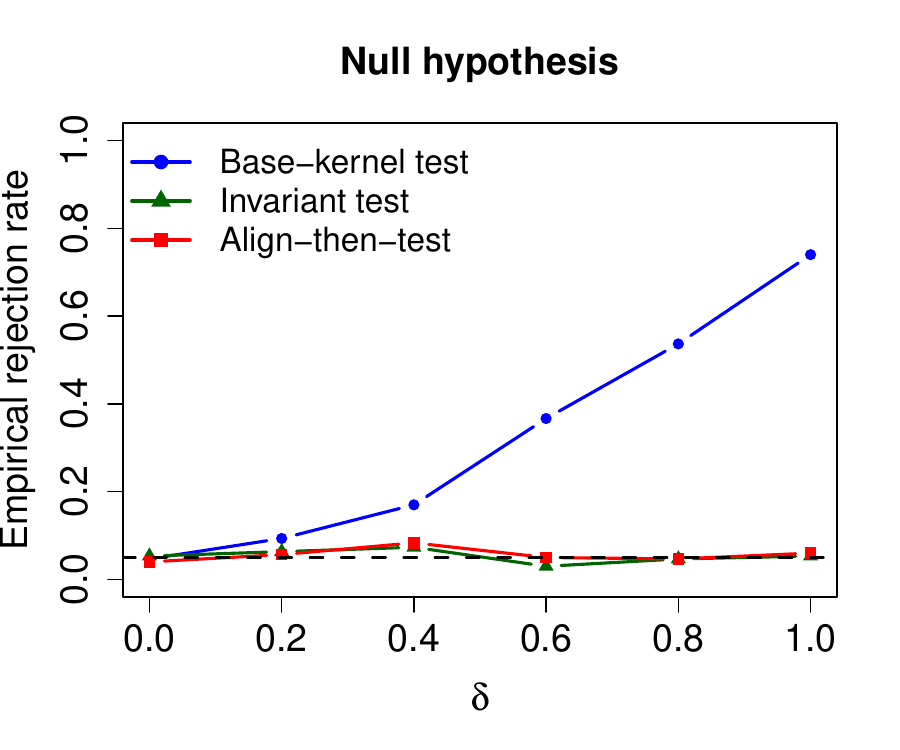}
    \hfill
    \includegraphics[width=0.49\textwidth]{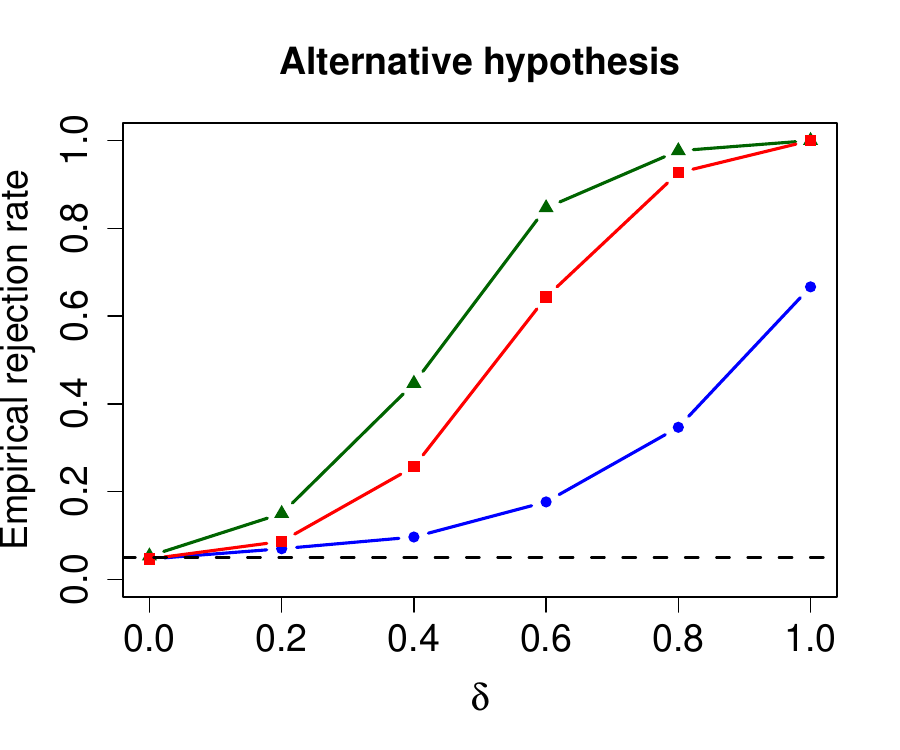}
    \caption{\textbf{Aperiodic case}. Empirical rejection rates of the MMD tests based on $k$, $k^\rho_\lambda$ and the align-then-test baseline with respect to $\delta$. The approximation budget of $k^\rho_\lambda$ is $S=16$. \textbf{Left.} $P \neq Q$ and $\Pi_*P = \Pi_*Q$. \textbf{Right.} $\Pi_*P\neq\Pi_*Q$.}
    \label{fig:exp2-aperiodic}
\end{figure}

\begin{figure}[H]
    \centering
    \includegraphics[width=0.48\textwidth]{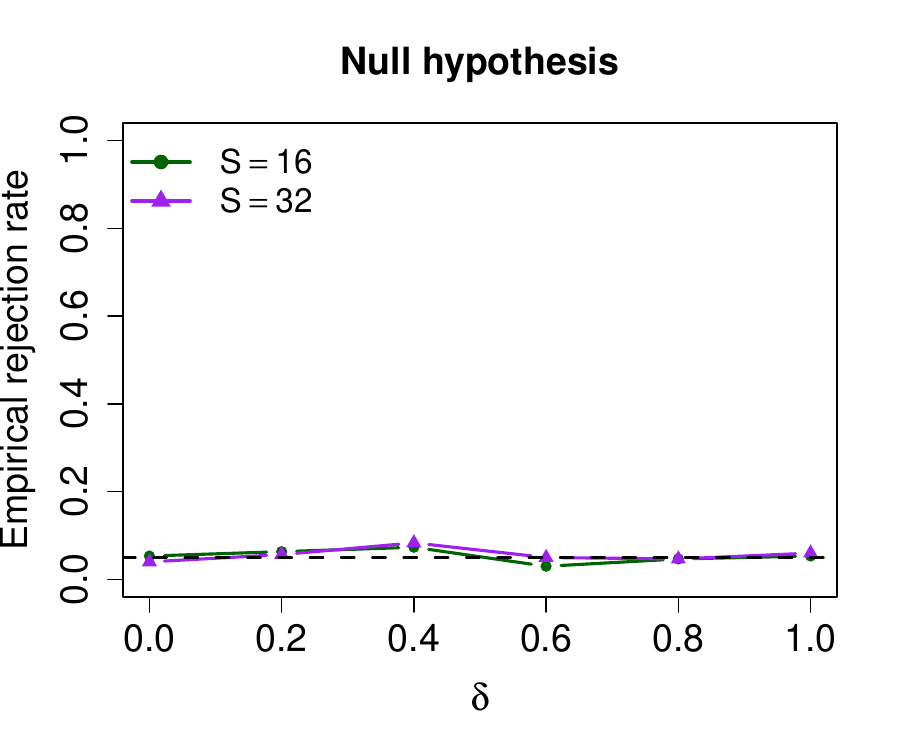}
    \hfill
    \includegraphics[width=0.48\textwidth]{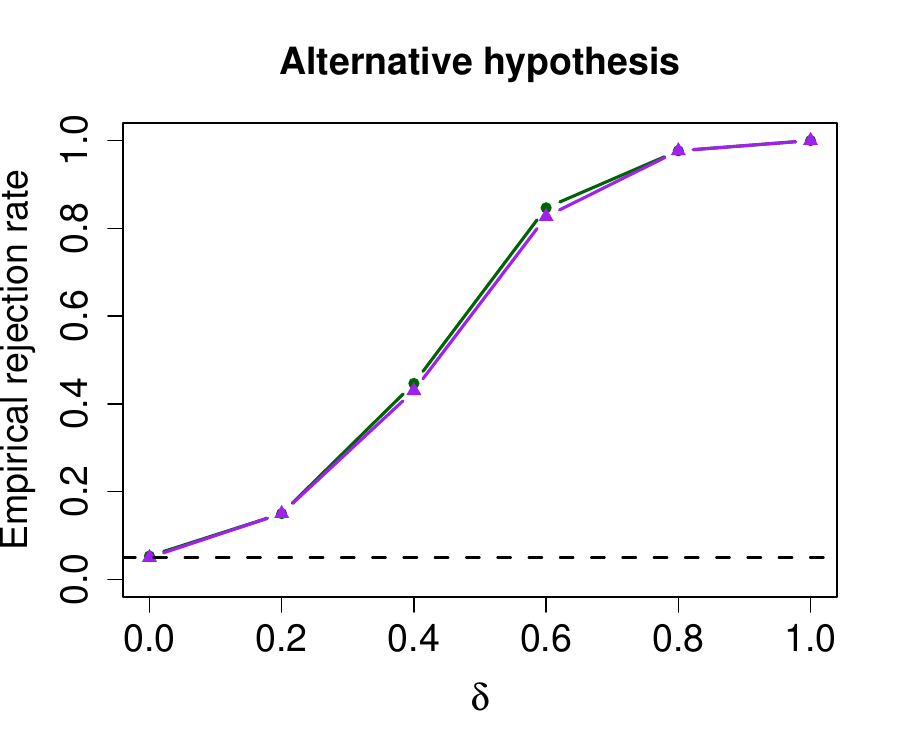}
    \caption{\textbf{Aperiodic case}. Empirical rejection rates of the MMD tests based on $k^\rho_\lambda$with two different approximation budgets. \textbf{Left.} $P \neq Q$ and $\Pi_*P = \Pi_*Q$. \textbf{Right.} $\Pi_*P\neq\Pi_*Q$.}
    \label{fig:exp2-aperiodic-S32}
\end{figure}

\medskip
Figure \ref{fig:exp2-aperiodic} shows the empirical rejection rates with respect to $\delta$ under the $\mathcal{H}_0$ and $\mathcal{H}_1$ scenarios. The results in the aperiodic setting are similar to those obtained in the periodic case. The invariant and align-then-test procedures remain well calibrated under $\mathcal{H}_0$ and their power increases with~$\delta$. The invariant test achieves the best performance, followed by the align-then-test procedure, and then the baseline test. As in the periodic setting of Section \ref{sec:simu_periodic}, we assess the sensitivity to the approximation of $k^\rho_\lambda$. Results obtained with $S = 32$ are presented in Figure \ref{fig:exp2-aperiodic-S32} and remain essentially unchanged.

\section{Real data application}
\label{sec:real_data}

We illustrate our invariant two-sample test on phonocardiogram (PCG) signals from the public \emph{PhysioNet/Computing in Cardiology Challenge 2016} dataset (available \href{https://physionet.org/content/challenge-2016}{here}). A phonocardiogram is an audio recording of the sounds produced by the heart over successive cardiac cycles. The training set is partitioned into six subsets (\texttt{training-a} to \texttt{training-f}). It contains a total of 3126 recordings whose durations range from a few seconds to more than one minute. The signals are collected under heterogeneous conditions, using different sensors and in diverse clinical and environmental settings. In this application, we restrict our analysis to the \texttt{training-e} subset.

\begin{figure}[H]
  \centering
  \includegraphics[width=0.9\linewidth]{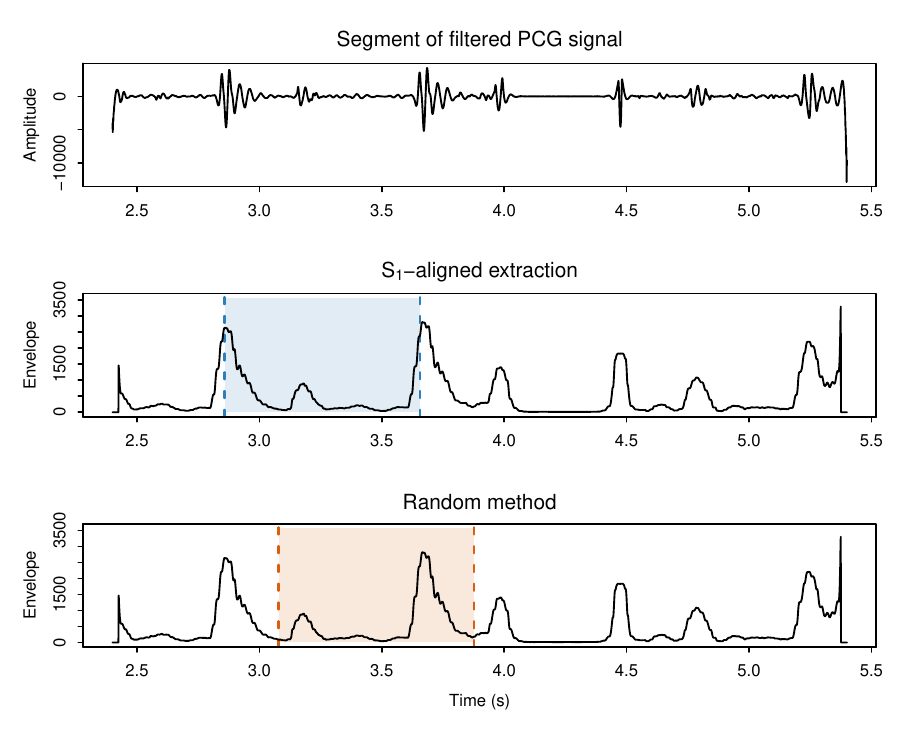}
  \caption{Illustration of the two cycle-extraction procedures on a single recording. In the $S_1$-aligned extraction, an interval $[s_1, s_1+\widehat{T})$ is extracted starting from a detected $s_1$ position. In the random extraction, an interval $[t_0, t_0+\widehat{T})$ is extracted from a randomly selected starting time $t_0$. Dashed vertical lines mark the beginning and end of the extracted interval.}
  \label{fig:extraction}
\end{figure}

\medskip
\textbf{Preprocessing and periodic representation.} All audio files are resampled to a common time grid with 1000 measurements per second. Then, a band-pass filter removes very slow trends and rapid fluctuations caused by noise. We compute a smoothed energy envelope of the signal using a sliding RMS (root mean square). Next, we compute the autocorrelation of the envelope. The characteristic cardiac period, denoted $\widehat{T}$, is the time lag that maximizes the autocorrelation within a plausible heart-rate range. The estimated period is used only to extract a fixed-length segment from each recording. This segment is treated as one period of a periodic signal. From each file, we construct a vector using one of the strategies below, illustrated in Figure \ref{fig:extraction}.

\medskip
\textbf{$S_1$-aligned extraction.} The first heart sound, denoted $\mathrm{S_1}$, typically appears as a prominent local peak in the energy envelope. 
We select one detected $S_1$ time point $s_1$ and extract the signal segment $[s_1,\,s_1+\widehat T)$, corresponding approximately to one cardiac cycle. The segment is then interpolated and resampled to obtain $p=128$ evenly spaced values. Finally, the resulting vector is standardized to have zero mean and unit variance.

\medskip
\textbf{Random-start extraction.} After estimating $\widehat{T}$, we draw a starting time $t_0$ uniformly among all positions where the segment $[t_0, t_0 + \widehat{T})$ lies entirely within the signal. This segment is then extracted, resampled to length $p=128$ and standardized. Unlike the $S_1$-aligned procedure, this method does not align the segment with a specific physiological landmark. 

\begin{figure}[H]
  \centering  \includegraphics[width=0.49\linewidth]{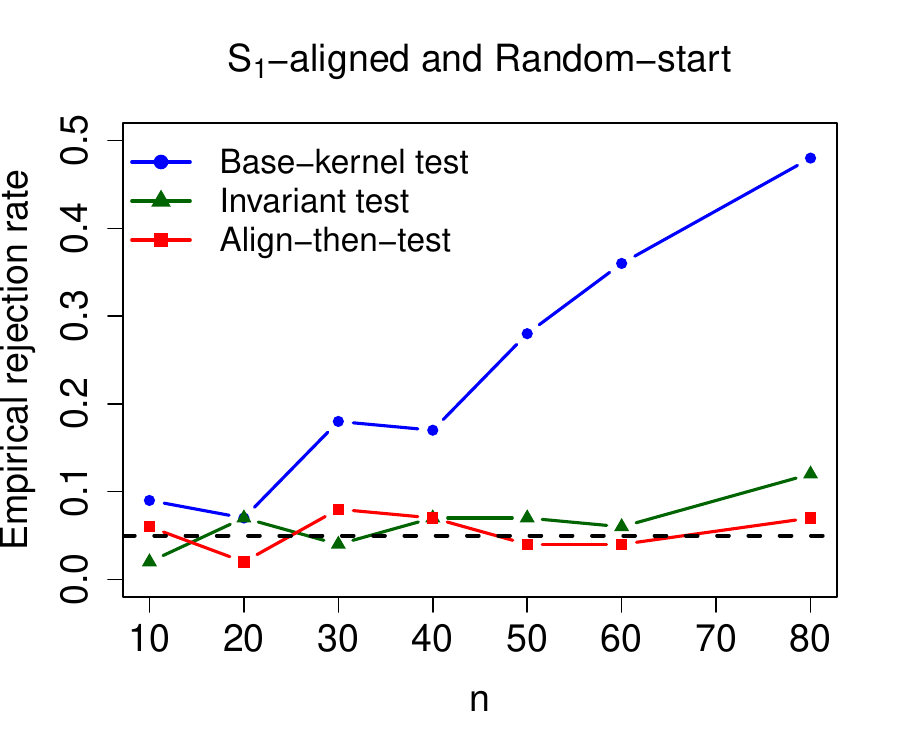}
  \caption{Empirical rejection rates of the MMD tests based on $k$, $k^\rho_\lambda$ and the align-then-test baseline with respect to $n$. The approximation budget of $k^\rho_\lambda$ is $S=16$.}
  \label{fig:realdata}
\end{figure}

\medskip
\textbf{Misalignment two-sample testing.} We consider sample sizes ranging from $n = 10$ to $n = 80$. For each value of $n$, we randomly select 
$2n$ distinct recordings from the \texttt{training-e} pool, without replacement. From each recording, a single signal is extracted using one of the two procedures described above. The sample $\mathbb{X}_n$ is obtained using the $S_1$-aligned extraction, while the sample $\mathbb{Y}_n$ is obtained using the random-start extraction. In this protocol, the distributions $P$ and $Q$ may differ because the extracted segments have different phase origins. Since all recordings originate from the same underlying pool, this difference is intended to reflect a temporal alignment mismatch. Consequently, after projection onto the quotient by circular translations, one expects $\Pi_* P$ and $\Pi_* Q$ to be approximately equal. Figure \ref{fig:realdata} shows the empirical rejection rates of each testing procedure with respect to the sample size $n$. We observe that the base-kernel test becomes more powerful as $n$ increases, reflecting its sensitivity to changes in the phase origin. In contrast, the invariant test and the align-then-test procedure are fairly insensitive to this phase mismatch. They exhibit similar behavior and low rejection rates. 

\begin{figure}[H]
  \centering
  \includegraphics[width=0.49\linewidth]{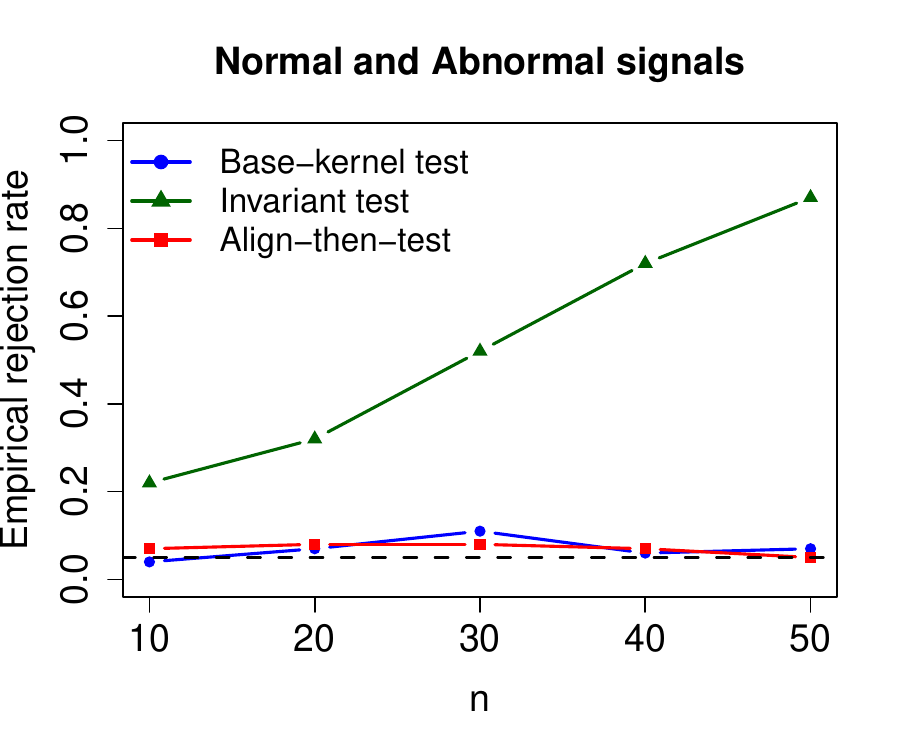}
  \caption{Empirical rejection rates of the MMD tests based on $k$, $k^\rho_\lambda$ and the align-then-test baseline with respect to $n$. The approximation budget of $k^\rho_\lambda$ is $S=16$.}
  \label{fig:realdata_normvsanorm}
\end{figure}

\medskip
\textbf{Two-Sample Testing for Normal and Abnormal Signals.} The PCG dataset provides, for each recording, a binary label normal/abnormal, indicating the absence or presence of cardiac anomalies. We use these labels to form two samples and extract the signals with the $S_1$-aligned procedure. Figure \ref{fig:realdata_normvsanorm} shows the empirical rejection rates with respect to $n$. The invariant test becomes more powerful as $n$ increases. This indicates that it detects distributional differences between the two samples. However, the power of the other two procedures remains close to the nominal level $\alpha=0.05$, independently of the value of $n$. This indicates that both procedures fail to capture the differences between the two sample distributions. As observed in the simulation study in Section \ref{sec:simulations}, translation variability can mask differences in signal morphology. In this case, the invariant-kernel approach can be substantially more powerful than the base-kernel test. In this real-data setting, the behavior of the align-then-test procedure suggests that the alignment step does not always remove phase variability. As a result, the subsequent kernel-based test may fail to detect the underlying shape differences.

\section{Proofs}
\subsection{Proof of Proposition~\ref{prop:S_nu}}
Let $\nu\in \mathcal{M}^\rho(G)$. By definition,
$$\int_{\mathcal{X}} \rho(x)\,\mathrm{d} (S_\nu P)(x)
=
\int_{\mathcal{X}}  \rho(x)\left(\int_G (\varphi_g)_*P(\mathrm{d}x) \, \mathrm{d} \nu(g) \right).$$
Recall that $\nu$ is $\sigma$-finite. By applying Tonelli's theorem, we obtain
\begin{align*}
\int_{\mathcal{X}} \rho(x)\,\mathrm{d}(S_\nu P)(x)
&=
\int_G \int_{\mathcal{X}} \rho(x)\,(\varphi_g)_*P(\mathrm{d} x) \,\mathrm{d}\nu(g)\\
&= \int_G \int_{\mathcal{X}} \rho(\varphi_g(x)) \,\mathrm{d}P(x) \, \mathrm{d}\nu(g)\\
&= \int_{\mathcal{X}} \left(\int_G \rho (\varphi_g(x) )\,\mathrm{d}\nu(g)\right)\,\mathrm{d}P(x).
\end{align*}
By definition of $\mathcal{M}^\rho(G)$, for all $x$ in $\mathcal{X}$, we have
$$
\int_G \rho (\varphi_g(x)) \,\mathrm{d} \nu (g)
\le
\sup_{z\in \mathcal{X}} \left\{ \int_G \rho (\varphi_g(z)) \,\mathrm{d} \nu (g) \right\} < \infty.$$
We now introduce
$$C_\nu := \sup_{z\in \mathcal{X}} \left\{ \int_G \rho (\varphi_g(z)) \,\mathrm{d} \nu (g) \right\}.$$
Then, 
$$\int_{\mathcal{X}} \rho(x)\,\mathrm{d} (S_\nu P)(x)
\le
\int_{\mathcal{X}} C_\nu\,\mathrm{d}P(x)
= C_\nu < \infty.$$
This shows that $S_\nu P \in \mathcal{M}^\rho(\mathcal{X})$.

\subsection{Proof of Proposition~\ref{prop:carac}}
Let $\mu\in \mathcal{M}^\rho(\mathcal{X})$ and let $\rho \mu$ be the measure defined for all  $A$ in $\mathcal B(\mathcal{X})$ by 
$$\rho\mu(A) := \int_A \rho(x)\,d\mu(x).$$
Since $\rho>0$ and $\displaystyle \int \rho\,\mathrm{d}|\mu|<\infty$, the measure $\rho\mu$ is well-defined and finite. By definition of the mean embedding $m^\rho (\mu)$ of $\mu$, associated with the kernel $k^\rho$, for all $x$ in $\mathcal{X}$ we have 
$$
m^\rho (\mu) (x) = \int_\mathcal{X} k^\rho (x,y)\,\mathrm{d}\mu(y) = \rho(x) \int_\mathcal{X} k(x,y) \rho(y) \, \mathrm{d}\mu(y) =  \rho(x) \int_\mathcal{X} k(x,y) \, \mathrm{d}(\rho\mu)(y).$$
In other words,
$$m^\rho (\mu) = \rho \times m (\rho \mu),$$
where $m (\rho \mu)$ denotes the mean embedding of $\rho \mu$ with the kernel $k$. Let now $\mu^1,\mu^2$ in $\mathcal{M}^\rho(\mathcal{X})$ such that $m^\rho (\mu^1) = m^\rho (\mu^2)$. Then, for all $x$ in $\mathcal{X}$,
$$\rho \times m (\rho \mu^1) = \rho \times m (\rho \mu^2).$$
Since $\rho > 0$, we obtain
$m (\rho \mu^1) = m (\rho \mu^2)$. Knowing that $k$ is characteristic on the space of finite signed measures $\mathcal{M}_f(\mathcal{X})$, it follows that
$$\rho\mu^1 = \rho\mu^2.$$
Let $\sigma := \mu^1-\mu^2$, we have $\rho\sigma=0$.
For $n\ge1$, set 
$$A_n:= \left\{ x\in \mathcal{X} \, | \, \rho(x) \geq \frac{1}{n} \right\}.$$ 
Then, $(A_n)_{n \geq 1}$ is an increasing nested sequence and
$\displaystyle \bigcup_{n\geq1} A_n = \mathcal{X}$.
For all $B$ in $\mathcal B(\mathcal{X})$, we have
$$\sigma(B\cap A_n)=\int_\mathcal{X}\frac{\mathbf 1_{B\cap A_n}}{\rho} \,\mathrm{d}(\rho\sigma)=0.$$
Hence, for all $n \geq 1$, $\sigma|_{A_n}=0$. Now define $C_1=A_1$ and $C_n=A_n \setminus A_{n-1}$ for $n \geq 2$. The sets $(C_n)_{n\ge1}$ are pairwise disjoint and $\displaystyle \bigcup_{n\ge1} C_n = \mathcal{X}$. Then, for all $B$ in $\mathcal B(\mathcal{X})$, 
$$\sigma(B)=\sum_{n\ge1}\sigma(B\cap C_n)=0,$$
since $B\cap C_n\subseteq A_n$ and $\sigma|_{A_n}=0$. Therefore $\sigma=0$ and $\mu^1=\mu^2$. This shows that the mean embedding $\mu \mapsto m^\rho (\mu)$ is injective on $\mathcal{M}^\rho(X)$, meaning that the kernel $k^\rho$ is characteristic on $\mathcal{M}^\rho(\mathcal{X})$.

\subsection{Proof of Proposition~\ref{prop:Av-kernel}}
We first show that $\Phi^\rho_\nu$ is well-defined. Recall that, $k$ and $\rho$ are continuous and bounded, the same applies to $k^\rho$. Then, by continuity of the group action, the mapping $g \mapsto \Phi^\rho(\varphi_g(x))$ is measurable. Denote by $\mathcal H^\rho_k$ the RKHS associated with $k^\rho$. Since $k^\rho$ is continuous and bounded and $\mathcal{X}$ is separable, then according to Corollary 4 of Section 1.5 in \cite{berlinet2011reproducing}, $\mathcal H^\rho_k$ is separable. Therefore, the mapping $g \mapsto \Phi^\rho(\varphi_g(x))$ is Bochner-measurable. In addition, for all $x \in \mathcal{X}$, we have
\begin{equation}
\label{eq:borne}
\|\Phi^\rho(\varphi_g(x))\|_{\mathcal H^\rho_k}
= \sqrt{k^\rho (\varphi_g(x),\varphi_g(x))}
= \sqrt{k(\varphi_g(x),\varphi_g(x))} \,\rho(\varphi_g(x))
\leq \sqrt{K}\,\rho(\varphi_g(x)).
\end{equation}
By definition of $\mathcal{M}^\rho(G)$,
$$\int_G \|\Phi^\rho(\varphi_g(x))\|_{\mathcal H^\rho_k}\,\mathrm{d} \nu(g)
\le \sqrt{K}\int_G \rho(\varphi_g(x))\,\mathrm{d} \nu (g) < \infty.$$
This shows that $g \mapsto \Phi^\rho(\varphi_g(x))$ is Bochner-integrable and that $\Phi^\rho_\nu$ is well-defined. Furthermore, for all $x,y\in \mathcal{X}$ we have
\begin{equation}
\label{eq:def_int}
\langle \Phi^\rho_\nu(x), \Phi^\rho_\nu(y)\rangle_{\mathcal H^\rho_k} 
= \left\langle \int_G  \Phi^\rho(\varphi_g(x))\, \mathrm{d}\nu(g),\, \int_G  \Phi^\rho(\varphi_h(y))\, \mathrm{d}\nu(h) \right\rangle_{\mathcal H^\rho_k}.
\end{equation}
Given that $\mathcal H^\rho_k$ is a separable Hilbert space and $\nu$ is $\sigma$-finite, we can apply Bochner-Fubini's theorem in~\eqref{eq:def_int} together with bilinearity of the inner product. Then, 
$$
\langle \Phi^\rho_\nu(x), \Phi^\rho_\nu(y)\rangle_{\mathcal H_k}
=
\int_G\!\int_G
\left\langle \Phi^\rho(\varphi_g(x)),  \Phi^\rho(\varphi_h(y)) \right\rangle_{\mathcal H^\rho_k} \mathrm{d}\nu(g)\,\mathrm{d}\nu(h).
$$
By the reproducing property of $\mathcal H^\rho_k$,
$$
\langle \Phi^\rho(\varphi_g(x)), \Phi^\rho(\varphi_h(y))\rangle_{\mathcal H^\rho_k}
= k^\rho(\varphi_g(x),\varphi_h(y)).
$$
Hence,
$$
\langle \Phi^\rho_\nu(x), \Phi^\rho_\nu(y)\rangle_{\mathcal H^\rho_k}
=
\int_G\! \int_G k^\rho(\varphi_g(x),\varphi_h(y))\mathrm{d}\nu(g)\, \mathrm{d}\nu(h)
= k^\rho_\nu(x,y),
$$
which proves the stated equality.

\subsection{Proof of Proposition~\ref{prop:Av-mean-emb}}
Using Equation \eqref{eq:borne}, we have
$$
\|\Phi^\rho(\varphi_g(x))\|_{\mathcal H^\rho_k}
\leq \sqrt{K}\,\rho(\varphi_g(x)).$$
Then, 
$$\int_\mathcal{X} \int_G \|\Phi^\rho (\varphi_g(x))\|_{\mathcal H^\rho_k}\, \mathrm{d}\nu (g) \,\mathrm{d}P(x)
\leq \sqrt{K} \int_\mathcal{X} \int_G \rho(\varphi_g(x))\, \mathrm{d}\nu (g) \, \mathrm{d}P(x) < \infty.$$
Indeed, $\nu \in M^\rho(G)$ and $P \in \mathcal P(X)$ and the same arguments of the proof of Proposition~\ref{prop:S_nu} hold. By definition of the mean embedding, 
\begin{align*}
m^\rho_\nu (P)
&= \int_\mathcal{X} \Phi^\rho_\nu(x)\,\mathrm{d}P(x) \\
&= \int_\mathcal{X} \int_G \Phi^\rho(\varphi_g(x))\,\mathrm{d}\nu(g)\,\mathrm{d}P(x)
\end{align*}
Using Bochner-Fubini's theorem yields to
\begin{align*}
m^\rho_\nu (P)
&= \int_G \left(\int_\mathcal{X} \Phi^\rho(\varphi_g(x))\,\mathrm{d}P(x)\right)\mathrm{d}\nu(g) \\
&= \int_G \left(\int_\mathcal{X} \Phi^\rho(u)\,\mathrm{d}((\varphi_g)_*P)(u)\right)\mathrm{d}\nu(g) \\
&= \int_\mathcal{X} \Phi^\rho(u)\,\mathrm{d}\left(\int_G (\varphi_g)_*P\,\mathrm{d}\nu(g)\right)(u) \\
&= \int_\mathcal{X} \Phi^\rho(u)\,\mathrm{d}(S_\nu P)(u) \\
&= m^\rho (S_\nu P).
\end{align*}
This concludes the proof.

\subsection{Proof of Theorem~\ref{prop:Haar-quotient}}
The Haar measure $\lambda$ is assumed to belong to $\mathcal{M}^\rho(G)$. Recall that, for all probability measure $P$ on $\mathcal{X}$ and all Borel set $A$, we have 
$$S_\lambda P(A) = \int_G P \left( \varphi_g^{-1}(A) \right) \mathrm{d} \lambda(g).$$

\medskip 
($\Rightarrow$)
Assume that $S_\lambda P = S_\lambda Q$.

\medskip
Let $\psi : \mathcal{X}/G \to \mathbb R$ be a bounded Borel function. First assume that $\psi\ge 0$. The general bounded case follows by writing
$\psi=\psi^+ - \psi^-$, where $\psi^+ , \psi^- \geq 0$ and by using integral linearity. We define
$$f : x \mapsto \psi(\Pi(x)) \rho(x).$$
Knowing that $S_\lambda P$ and $S_\lambda Q$ belong to $\mathcal{M}^\rho(\mathcal{X})$, the integrals of $f$ with respect to these measures are finite and satisfy
$$\int_\mathcal{X} f \, \mathrm{d}(S_\lambda P) = \int_\mathcal{X} f \, \mathrm{d}(S_\lambda Q).$$
In addition, 
$$\int_\mathcal{X} f(x)\,\mathrm{d}(S_\lambda P)(x)
=
\int_G\int_\mathcal{X} f(x)\,\mathrm{d}(\varphi_g)_*P(x)\, \mathrm{d}\lambda(g)
=
\int_G\int_\mathcal{X} f(\varphi_g(x))\,\mathrm{d}P(x)\, \mathrm{d}\lambda(g),$$
where we used the definition of $S_\lambda P$ and a change of variables.
Now,
$$f(\varphi_g(x)) = \psi(\Pi(\varphi_g(x)))\,\rho(\varphi_g(x)) = \psi(\Pi(x))\,\rho(\varphi_g(x)),$$
since $\Pi\circ \varphi_g = \Pi$.
Therefore,
$$\int_\mathcal{X} f\,\mathrm{d}S_\lambda P =
\int_\mathcal{X} \psi(\Pi(x))
\left( \int_G \rho(\varphi_g(x))\,\mathrm{d}\lambda(g) \right)\, \mathrm{d}P(x).$$
Define
$$w : x \mapsto \int_G \rho(\varphi_g(x))\, \mathrm{d}\lambda(g).$$
Since $\lambda \in \mathcal{M}^\rho(G)$, the function $w$ is finite and bounded on $\mathcal{X}$, and Haar invariance implies that $w \circ \varphi_h =w$, for all $h$ in $G$. Thus, there exists a measurable function $\tilde w : \mathcal{X}/G \to (0,\infty)$ such that $w = \tilde w \circ \Pi$.
We obtain
$$\int_\mathcal{X} f\,\mathrm{d}(S_\lambda P)
=
\int_\mathcal{X} \psi(\Pi(x))\,\tilde w(\Pi(x))\,\mathrm{d}P(x)
=
\int_{\mathcal{X}/G} \psi(y)\,\tilde w(y)\,\mathrm{d}(\Pi_*P)(y).$$
Similarly,
$$\int_\mathcal{X} f\,\mathrm{d}(S_\lambda Q)
=
\int_{\mathcal{X}/G} \psi(y)\,\tilde w(y)\, \mathrm{d}(\Pi_*Q)(y).$$
Since $S_\lambda P = S_\lambda Q$, the two integrals are equal for all bounded Borel $\psi$. Hence,
$$
\tilde w\,\Pi_*P = \tilde w\,\Pi_*Q,$$
as measures on $\mathcal{X}/G$. Let $\eta:=\Pi_*P-\Pi_*Q$, we have $\tilde w\,\eta=0$. For $n\geq 1$, set 
$$B_n:= \left\{y \in \mathcal{X}/G \; | \; \widetilde w(y) \geq \frac{1}{n} \right\},$$
Then, $(B_n)_{n \geq 1}$ is an increasing nested sequence and
$\displaystyle \bigcup_{n\geq1} B_n = \mathcal{X}/G$. For all $B$ in $\mathcal B(\mathcal{X}/G)$, 
$$
\eta(B\cap B_n)
=\int_{\mathcal{X}/G} \frac{\mathbf 1_{B\cap B_n}}{\tilde w}\, \mathrm{d}(\tilde w\,\eta)=0,$$
since $\mathbf 1_{B\cap B_n}/\tilde w \leq n$.
Hence $\eta|_{B_n}=0$ for all $n \geq 1$. Therefore, since $\eta$ is a finite signed measure, $\eta=0$. In other words, $$\Pi_*P=\Pi_*Q.$$

\medskip\noindent
($\Leftarrow$)
Assume that $\Pi_*P = \Pi_*Q$.

\medskip
For a Borel set $A$ in $\mathcal{B}(\mathcal{X})$, define
$$F_A : x \mapsto \int_G \mathbf 1_A(\varphi_g(x))\, \mathrm{d}\lambda(g).$$
By Haar invariance, $F_A \circ \varphi_h  = F_A$ for all $h\in G$. Then, $F_A$ is constant on orbits.
Hence, there exists a measurable function $f_A : \mathcal{X}/G \to [0,\infty]$ such that $F_A = f_A\circ \Pi$. By definition of $S_\lambda P$ and Tonelli's theorem,
$$S_\lambda P(A)
=
\int_G P(\varphi_g^{-1}(A))\, \mathrm{d}\lambda(g)
=
\int_\mathcal{X} F_A(x)\, \mathrm{d}P(x)
=
\int_{\mathcal{X}/G} f_A(y)\,\mathrm{d}(\Pi_*P)(y).$$
Similarly,
$$S_\lambda Q(A) =
\int_{\mathcal{X}/G} f_A(y)\, \mathrm{d}(\Pi_*Q)(y).$$
Since $\Pi_*P = \Pi_*Q$, we have $S_\lambda P(A) = S_\lambda Q(A)$ for all $A$ in $\mathcal{B}(\mathcal{X})$. Meaning that, 
$$S_\lambda P = S_\lambda Q.$$

\medskip\noindent
This concludes the proof of Theorem \ref{prop:Haar-quotient}.

\section*{Acknowledgments}

The authors acknowledge the support of the French National Research Agency (ANR) through the project ANR-24-CE40-2439 (FUNMathStat). They also thank Magalie Fromont and Nicolas Klutchnikoff for valuable discussions.

\bibliographystyle{unsrt}
\bibliography{Refs}
\end{document}